\documentclass[11pt]{article}%
\usepackage{geometry}                
\usepackage{tikz} 
\usetikzlibrary{positioning}
\usepackage{graphicx}
\usepackage{amssymb}
\usepackage{epstopdf}
\usepackage{subfigure}  
\usepackage[sans]{dsfont}
\usepackage[english]{babel}
\usepackage{latexsym}
\usepackage{mathrsfs}
\usepackage{graphicx}
\usepackage{color}
\usepackage{float}
\usepackage{mathtools}
\usepackage{amsmath}
\usepackage{amsfonts}
\numberwithin{equation}{section}
\usepackage{enumerate}
\usepackage{amsthm}
\usepackage[title]{appendix}

\usepackage[bookmarks=true,colorlinks=true,linkcolor={blue},urlcolor={blue}, citecolor={blue},pdfstartview={XYZ null null 1.22}]{hyperref}%

\newcommand{\R}{\mathbb R}

\def\be#1\ee{\begin{equation}#1\end{equation}}

\theoremstyle{definition}

\newtheorem{remark}{Remark}

\def\be{\begin{equation}}
	\def\ee{\end{equation}}
\def\bea{\begin{eqnarray}}
	\def\eea{\end{eqnarray}}

\newenvironment{equations}{\equation\aligned}{\endaligned\endequation}
\title{Augmented data and neural networks for robust epidemic forecasting: application to COVID-19 in Italy}
\date{} 

\author{ G. Dimarco,\thanks{Department of Mathematics and Computer Science \& Center for Modeling, Computing and Statistics (CMCS), University of Ferrara, via Machiavelli 30, 44121 Ferrara, ITALY. (giacomo.dimarco@unife.it)} \and
	F. Ferrarese, \thanks{Department of Mathematics and Computer Science \& Center for Modeling, Computing and Statistics (CMCS), University of Ferrara, via Machiavelli 30, 44121 Ferrara, ITALY. (federica.ferrarese@unife.it)} \and
	L. Pareschi \thanks{Department of Mathematics and Computer Science \& Center for Modeling, Computing and Statistics (CMCS), University of Ferrara, via Machiavelli 30, 44121 Ferrara, ITALY. (lorenzo.pareschi@unife.it) \&
		Maxwell Institute and Department of Mathematics
		School of Mathematical and Computer Sciences
		Heriot-Watt University, Edinburgh, UK (L.Pareschi@hw.ac.uk)
	}
}

\begin{document}
	\maketitle
	\section*{Abstract}
	In this work, we propose a data augmentation strategy aimed at improving the training phase of neural networks and, consequently, the accuracy of their predictions. Our approach relies on generating synthetic data through a suitable compartmental model combined with the incorporation of uncertainty. The available data are then used to calibrate the model, which is further integrated with deep learning techniques to produce additional synthetic data for training.
	The results show that neural networks trained on these augmented datasets exhibit significantly improved predictive performance. We focus in particular on two different neural network architectures: Physics-Informed Neural Networks (PINNs) and Nonlinear Autoregressive (NAR) models. 
	The NAR approach proves especially effective for short-term forecasting, providing accurate quantitative estimates by directly learning the dynamics from data and avoiding the additional computational cost of embedding physical constraints into the training. In contrast, PINNs yield less accurate quantitative predictions but capture the qualitative long-term behavior of the system, making them more suitable for exploring broader dynamical trends. 
	Numerical simulations of the second phase of the COVID-19 pandemic in the Lombardy region (Italy) validate the effectiveness of the proposed approach.
	\begin{center}
		{\bf Keywords:} {Mathematical epidemiology, Physics-Informed Neural Networks, Non-linear Autoregressive Deep learning, COVID 19 data}
	\end{center}

\section{Introduction}\label{sec:intro} 
Throughout history, pandemics have had a profound impact on global health, economies, and everyday life \cite{chowell2006transmission, dimarco2020wealth,lunelli2009epidemic,toscani2022multi}. From the Spanish flu in 1918 to more recent outbreaks such as the H1N1 influenza in 2009 and the COVID-19 pandemic, the repeated emergence of infectious diseases has emphasized the urgent need for timely and accurate response strategies \cite{albi2021control, albi2021modelling, bolzoni2019optimal,flaxman2020estimating,lee2010optimal} .
In this context, mathematical modeling plays a crucial role in predicting disease transmission, assessing the effectiveness of intervention measures, and informing public health decisions. Among the most commonly used approaches are traditional compartmental models, such as the Susceptible-Infected-Recovered (SIR) framework, valued for their simplicity and interpretability \cite{capasso1978,hethcote2000,kermack1927}. These models classify the population into key compartments: Susceptible individuals who are at risk of infection, Infected individuals who can spread the disease, and Recovered individuals who have either recovered or died and are no longer infectious. The basic SIR model can be extended to incorporate additional compartments \cite{chisholm2018implications, anastassopoulou2020data}. For example, the SIAR model introduces an asymptomatic class, especially relevant for diseases like COVID-19, where asymptomatic individuals significantly contribute to transmission \cite{leung2018infector, mizumoto2020estimating}. Despite their simplicity, such models often fall short in capturing the full complexity and variability of real-world epidemic dynamics, especially when dealing with incomplete or uncertain data. In the early stages of an outbreak, underreporting of infections is common, making it essential to incorporate uncertainty in model parameters or initial conditions to achieve more realistic scenarios \cite{albi2022kinetic, bertaglia2021spatial, capaldi2012parameter}. Additionally, the simplifying assumptions underlying traditional models can limit their ability to reflect population heterogeneity and the evolving nature of outbreaks. To address these limitations, more sophisticated models have been proposed that allow for time and state dependent transmission rates \cite{dimarco2021kinetic, korobeinikov2005non, zanella2021data}.  These enhancements enable a better representation of intervention measures—such as lockdowns which, in the case of COVID-19, played a significant role in reducing transmission. Moreover, accounting for individual heterogeneity is essential to capture the varying behaviors observed among individuals from different groups, such as distinct age classes \cite{castillo1989epidemiological, franceschetti2008threshold, voinsky2020effects}. 

In recent years, the use of neural networks for epidemic prediction has emerged as a promising complement to traditional modeling approaches \cite{han2024approaching, lu2021bifidelity, bertaglia2022asymptotic, Goodfellow-et-al-2016, millevoi2024physics}. In particular, neural networks for time series forecasting have shown strong potential. While purely data-driven models such as feed-forward networks can be useful for interpolation, they often lack physical consistency \cite{mcinerney2024statistical, pashaei2021training, chen2000new}. In contrast, Physics-Informed Neural Networks (PINNs) incorporate the governing equations into the training  \cite{de2024numerical, berkhahn2022physics, han2024approaching}, leading to more realistic solutions. However, the process involves automatic differentiation and other intrinsic computations which may led to very computationally expensive simulations. Alternative architectures, such as Recurrent Neural Networks (RNNs), offer a different approach by directly learning temporal dependencies without explicitly embedding the physical model into the training process \cite{mienye2024recurrent, li2020recurrent, amendolara2023lstm}. A notable example is the Nonlinear Autoregressive (NAR) network, which, while maintaining a feed-forward structure, effectively captures time dynamics by using past observations to predict future values \cite{sarkar2019comparative}. These networks offer a data-driven and computationally efficient alternative to PINNs for short-term forecasting, delivering a quantitative description of disease progression that can be particularly valuable for monitoring and managing the pandemic within hospitals. In contrast, PINNs provide a more qualitative representation of the disease dynamics, which can be especially useful for investigating strategies aimed at mitigating epidemic peaks during a pandemic.

In this work, we consider a specific compartmental model that has been shown to effectively capture the time evolution of epidemic spread in the presence of uncertain parameters \cite{zanella2021data}. This model extends a SIAR-type framework by incorporating age-structured dynamics to better reflect the pandemic’s impact across different demographic groups, as well as the influence of lockdown measures. Starting from this depicted model, we use real-world data to calibrate the model parameters under uncertainty \cite{capaldi2012parameter, chen2000new}.

We then introduce two different neural network architectures, namely Physics-Informed Neural Networks (PINNs) and Nonlinear Autoregressive (NAR) networks, and investigate how data augmentation strategies can enhance their predictive performance. In our previous work \cite{awais2025data}, we explored such strategies in a simplified context of deterministic epidemic models by augmenting training datasets with synthetic data generated from model simulations. Here, we extend this analysis to models containing uncertainty.

Through a series of numerical experiments, we demonstrate that NAR networks can accurately capture epidemic dynamics in both interpolation and extrapolation tasks. In short-term forecasting, they outperform PINNs, particularly when trained on augmented datasets. We also highlight the computational advantages of NAR networks: unlike PINNs, which require evaluating the underlying differential equations during training, NAR networks rely solely on data, leading to faster training and lower computational costs. In contrast, for long-term forecasting, PINNs provide more reliable predictions, particularly in capturing the epidemic peaks.

The rest of the paper is organized as follows:
In Section \ref{sec:problem_setting}, we introduce the mathematical models, starting from the classical SIR, we introduce a recent compartmental model which takes into account social behavior, age-structure and uncertainty \cite{zanella2021data}.
In Section \ref{sec:params_est}, we detail the parameter estimation procedure under uncertainty using real world data which permits to match the model evolution with the COVID-19 epidemic spread. Section \ref{sec:neural_networks} presents the neural network architectures used for epidemic prediction: PINNs and NAR networks trained with data and models outcome.
In Section \ref{sec:numerical_exp}, we report different numerical experiments to evaluate the networks performance, making a comparison between PINNs and NARs in short and long term forecasting.
Finally, Section \ref{sec:conclusion} outlines possible future research directions.

\section{Model setting: compartments, social behavior, age structure and uncertainty} \label{sec:problem_setting}
In our analysis, we will consider a suitable extension of the classical SIR model \cite{capasso1978, kermack1927}
\begin{equation}\label{eq:deterministicSIR}
	\begin{split}
		&\frac{dS(t)}{dt} = -\beta \frac{S(t) I(t)}{N},\\
		&\frac{dI(t)}{dt}  = \beta \frac{S(t) I(t)}{N} -\gamma I(t),\\
		&\frac{dR(t)}{dt} = \gamma I(t),
	\end{split}
\end{equation}  
which typically describes the spread of an infectious disease in a population of size $N$ by partitioning it into three compartments: Susceptible (S), Infected (I), and Recovered (R). While this model has been widely adopted in the past due to its simplicity, it has been proved to fall short in capturing complex epidemic dynamics \cite{Gatto,albi2022kinetic}, as it assumes, among other simplifying assumptions, the transmission and recovery rates, namely the parameters $\beta > 0$ and $\gamma > 0$, to be constant. Some significant modifications were recently proposed in \cite{dimarco2021kinetic}, where the authors, starting from a microscopic interaction dynamics, derived a compartmental model that accounts for the role of social contacts among individuals in the spread of an epidemic. By adding an additional variable characterizing the number of social contacts $x \ge0$ among individuals and by denoting with $f_S(x,t)$, $f_I(x,t)$ and $f_R(x,t)$, the distributions at time $t > 0$ of the number of social contacts of the population of susceptible, infected and recovered individuals, one can fix, upon renormalization, the total distribution of social contacts to be a probability density for all times $t \ge 0$
\[
\int_{\mathbb{R}_+} f(x,t)\,dx = 1.
\]
The model then follows combining the epidemic process with the contact dynamics. 
This gives the system 
\begin{equations}\label{sir-gamma}
	\frac{\partial f_S(x,t)}{\partial t} &= -K(f_S,f_I)(x,t) +   Q_S(f_S)(x,t),
	\\
	\frac{\partial f_I(x,t)}{\partial t} &= K(f_S,f_I)(x,t)  - \gamma f_I(x,t) + Q_I(f_I)(x,t),
	\\
	\frac{\partial f_R(x,t)}{\partial t} &= \gamma f_I(x,t) +  Q_R(f_R)(x,t),
\end{equations}
where the transmission of the infection is governed by the function
$K(f_S,f_I)$, the local incidence rate, expressed by
\be\label{inci}
{K(f_S,f_I)(x, t) = f_S(x,t) \int_{\R^+} \kappa(x,y)f_I(y,t) \,dy},
\ee
where $\kappa(x,y)$ is called the contact function, a nonnegative function growing with respect to the number of contacts $x$ and $y$ of the populations of susceptible and infected. A leading example reads
\[
\kappa(x,y) = \beta\, x^\alpha y^\alpha,
\]
where $\alpha, \beta$ are positive constants, that is by taking the incidence rate dependent on the product of the number of contacts of susceptible and infected people. The operators $Q_J$, $J\in \{S,I,R\}$ 
are integral operators that modify the distribution of contacts $f_J(x,t)$, $J\in \{S,I,R\}$ through repeated interactions among individuals \cite{dimarco2021kinetic}. Now, by integrating over the number of social contact and under suitable hypothesis 
on the operators $Q_J$, $J\in \{S,I,R\}$ characterizing the distribution of social contacts at equilibrium, one can observe that the evolution of the mass fractions obeys to a SIR-type model 
\begin{equation}\label{eq:deterministicSIR_H}
	\begin{split}
		&\frac{dS(t)}{dt} = -\beta S(t) I(t) H(I(t)),\\
		&\frac{dI(t)}{dt}  = \beta S(t) I(t)  H(I(t)) -\gamma I(t),\\
		&\frac{dR(t)}{dt} = \gamma I(t).
	\end{split}
\end{equation}  
Here, for simplicity, the population size $N$ is absorbed into the coefficient $\beta$, while the function $H(I(t))$ represents the average behavior induced by microscopic interactions \eqref{sir-gamma}. It denotes a macroscopic incidence rate that captures time-dependent modifications to the transmission dynamics, reflecting behavioral responses and public health interventions such as lockdowns. The inclusion of a state-dependent transmission rate enables the model to account for a broader range of epidemic scenarios as shown in \cite{dimarco2021kinetic}.

Moving forward into the modeling, we consider an additional extension to the framework \eqref{eq:deterministicSIR_H}, capable of addressing more complex scenarios, as shown in \cite{zanella2021data}. That is we include an additional compartment, 
$A$, to account for asymptomatic individuals, and more important we incorporate uncertainty into the model parameters. These extensions are particularly relevant for diseases such as COVID-19, where asymptomatic transmission it has been shown to play a crucial role \cite{albi2022kinetic}. Moreover, incorporating uncertainty provides a more realistic representation of data limitations, especially during the early stages of the pandemic, when the true number of infections was often significantly underreported. This permits to enhance the forecasting capabilities of the Neural Networks as shown later in the article. The resulting model is given by
\begin{equation}\label{eq:social_SIAR}
	\begin{split}
		&\frac{\partial S(t,z)}{\partial t} = -\Lambda(t,z),\\
		&\frac{\partial I(t,z)}{\partial t} = \xi(z) \Lambda(t,z) - \gamma_I(z) I(t,z),\\
		&\frac{\partial A(t,z)}{\partial t} = (1-\xi(z)) \Lambda(t,z) - \gamma_A(z) A(t,z),\\
		&\frac{\partial R(t,z)}{\partial t} =\gamma_I(z) I(t,z) + \gamma_A(z) A(t,z),
	\end{split}
\end{equation}
being 
\begin{equation}\label{eq:Lambda}
	\Lambda(t,z) =	\beta(z) S(t,z) H_S(I(t,z)) (H_I(I(t,z)) I(t,z) + H_A(I(t,z)) A(t,z)),
\end{equation}
where $\beta(z),\gamma_I(z), \gamma_A(z)>0$ are the transmission and recovery rates, and the macroscopic incidence rates $H_J(\cdot)$ for $J\in\{S,I,A\}$ are given by
\begin{equation}\label{eq:H}
	H_S(r) = \frac{\mu(z)}{\sqrt{1+\nu(z) r}},\qquad H_A(r) = H_S(r), \qquad H_I(r) = k H_S(r),
\end{equation}
with $\mu(z),\nu(z)>0$. These incidence rates model how, in response to the spread of the disease, both the susceptible and asymptomatic populations tend to reduce their average number of daily social contacts in a similar manner, while the contact rate of the infected population is further reduced by an additional factor $k\in[0,1]$. The parameter $\xi(z)$ represents the percentage of asymptomatic individuals in the population and takes values in the interval $[0,1]$.  The parameter $z$ represents the uncertainty and it is distributed according to a given distribution $p(z)$ as it will be specified later on.  

A further extension of the model in \eqref{eq:social_SIAR} consists in incorporating an age structure. In this case, the population is divided into subclasses corresponding to different age groups. As before, we consider four compartments: Susceptible (S), Infected (I), Asymptomatic (A), and Recovered (R). Each variable now depends not only on time and the uncertainty parameter $z$ but also on the age variable $x\in \mathcal{A} = (0,100)$ representing the age class. Incorporating an age structure introduces additional heterogeneity into the system, which is essential for accurately capturing the dynamics of an epidemic. Individuals belonging to different age groups typically display distinct contact patterns and social behaviors, resulting in varying transmission dynamics across the population. Under these assumptions, the model can be written as follows
\begin{equation}\label{eq:social_SIAR_ages}
	\begin{split}
		&\frac{\partial S(x,t,z)}{\partial t} = -\Lambda(x,t,z),\\
		&\frac{\partial I(x,t,z)}{\partial t} = \xi(x,z) \Lambda(x,t,z) - \gamma_I(x,z) I(x,t,z),\\
		&\frac{\partial A(x,t,z)}{\partial t} = (1-\xi(x,z)) \Lambda(x,t,z) - \gamma_A(x,z) A(x,t,z),\\
		&\frac{\partial R(x,t,z)}{\partial t} =\gamma_I(x,z) I(x,t,z) + \gamma_A(x,z) A(x,t,z),
	\end{split}
\end{equation}
being
\begin{equation}\label{eq:Lambda_ages}
	\begin{split}
		\Lambda(t,z,x) = &	\beta(x,z) S(x,t,z) H_S(I(x,t,z)) \\ &\qquad \qquad \int_{\mathcal{A}} (H_I(I(y,t,z)) I(y,t,z) + H_A(I(y,t,z)) A(y,t,z)) dy,
	\end{split}
\end{equation}
where  $\beta(x,z),\gamma_I(x,z), \gamma_A(x,z)>0$ are the transmission and recovery rates, now dependent also on the age class $x$, and the macroscopic incidence rates $H_J(\cdot)$ for $J\in\{S,I,A\}$ are given by
\begin{equation}\label{eq:H_ages}
	H_S(r) = \frac{\mu(x,z)}{\sqrt{1+\nu(x,z) r}},\qquad H_A(r) = H_S(r), \qquad H_I(r) = k H_S(r),
\end{equation}
where the age-dependent parameters $\mu(x,z),\nu(x,z)$ are supposed to be positive, for $k\in[0,1]$.

\section{Parameters estimation}\label{sec:params_est}
The goal of this section is to illustrate the calibration of the model parameters from available data. The procedure is described for the age-structured social SIAR model \eqref{eq:social_SIAR_ages}, but similar results can be obtained for the social SIAR model \eqref{eq:social_SIAR} without age structure by integrating over the age classes.
We focus on the second wave of COVID-19 in Italy, and in particular in the Lombardy region. Following the methodology outlined in \cite{zanella2021data}, we first estimate the recovery rates $\gamma_I$ and $\gamma_A$. These rates are derived from the Time of Viral Clearance (TVC), which measures the duration between the first positive test and recovery/death. As reported in \cite{zanella2021data}, the recovery times follow a beta distribution, allowing the parameters to be defined as follows
\begin{equation}\label{eq:gamma}
	\frac{1}{\gamma_I(x,z_1,z_2)} = \begin{cases}
		&h_{1,1} + h_{1,2}z_1, \qquad \text{if } x\leq 50,\\
		&h_{2,1} + h_{2,2}z_2,\qquad \text{if } x> 50,
	\end{cases}\qquad \gamma_A(x,z_1,z_2) = 2\gamma_I(x,z_1,z_2),
\end{equation}
being $h_{1,1} = h_{2,1} = 5$, $h_{1,2} = 32$, $h_{2,2} = 40$, and with $z_1\sim B(\alpha_1,\beta_1)$,
$z_2\sim B(\alpha_2,\beta_2)$ beta distributions with parameters $\alpha_1 = 2.1$, $\beta_1 = 5.1$, and $\alpha_2 = 1.8$, $\beta_2 = 3.9$. 
Next, the transmission rate $\beta(x,z)$ and the parameter $\xi(x,z)$ which represents the initial proportion of asymptomatic individuals, are estimated using data over the time interval  $[t_0,t_L]$, being $t_0$ the initial time, here supposed to be October 8th 2020, and $t_L$ the final time, here supposed to be October 20th 2020. In this period of time, there were no specific restrictions to the mobility and life on individuals.
To count for uncertainty, we generate samples $\{\gamma_I\}_{m=1}^M$ from $\gamma_I$, and  $\{\gamma_A\}_{m=1}^M$ from $\gamma_A$. These samples are constructed using a collocation approach based on Gauss-Jacobi polynomials with $M = 5$ nodes. For any $m=1,\ldots,M$, we solve the following optimization problem 
\begin{equation}\label{eq:min_beta}
	\min_{\beta_{m}(x),\xi_{m}(x)} \int_{\mathcal{A}} \mathcal{J}(I(y,t,z_m),\hat{I}(y,t),R(y,t,z_m),\hat{R}(y,t)) d y,  \qquad \text{for } t\in[t_0,t_L], 
\end{equation}
being 
\begin{equation}\label{eq:functional}
	\mathcal{J}(\cdot)  = p \Vert I(y,t,z_m) - \hat{I}(y,t) \Vert_{2} +(1-p) \Vert R(y,t,z_m) - \hat{R}(y,t) \Vert_{2},
\end{equation}
for $p\in[0,1]$, where $I(\cdot), R(\cdot)$ are the numerical solution to the underline model assuming $H_S =1$, meaning that no specific restrictions are set during this period of time, while $\hat{I}(\cdot), \hat{R}(\cdot)$ are the available data.  

\begin{figure}[h] 
	\centering
	\includegraphics[width=0.4\linewidth]{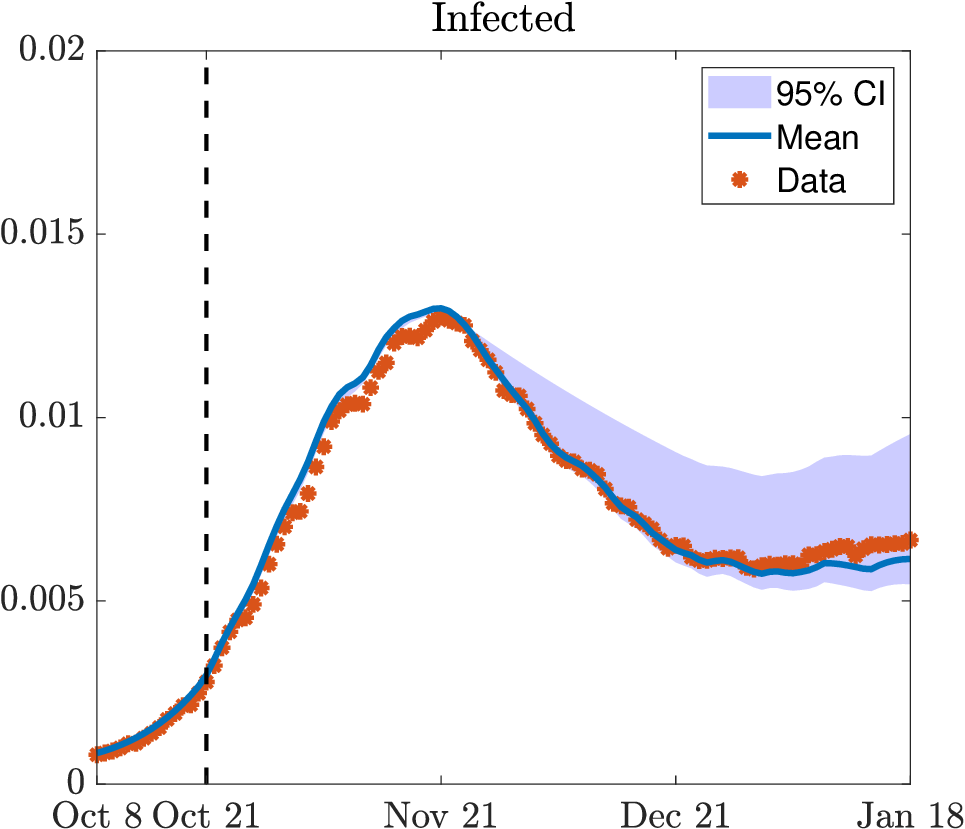} 
	\caption{Dynamics of the infected population obtained by solving the calibrated social-SIAR model \eqref{eq:social_SIAR} and compared with experimental data. The figure shows the mean epidemic trajectory with the 95\% confidence interval (shaded area), alongside the observed data. The black dashed line separates the two epidemic phases. }
	\label{fig:infected_noAges}
\end{figure}

\begin{figure}[h] 
	\centering
	\includegraphics[width=0.328\linewidth]{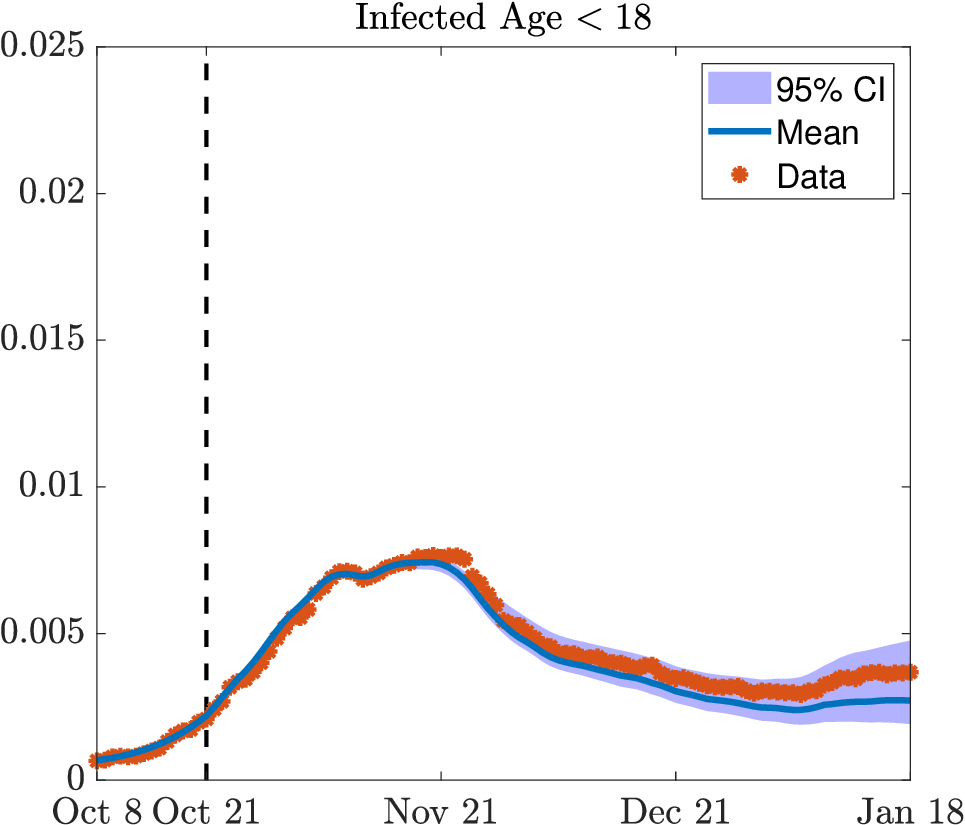 }
	\includegraphics[width=0.328\linewidth]{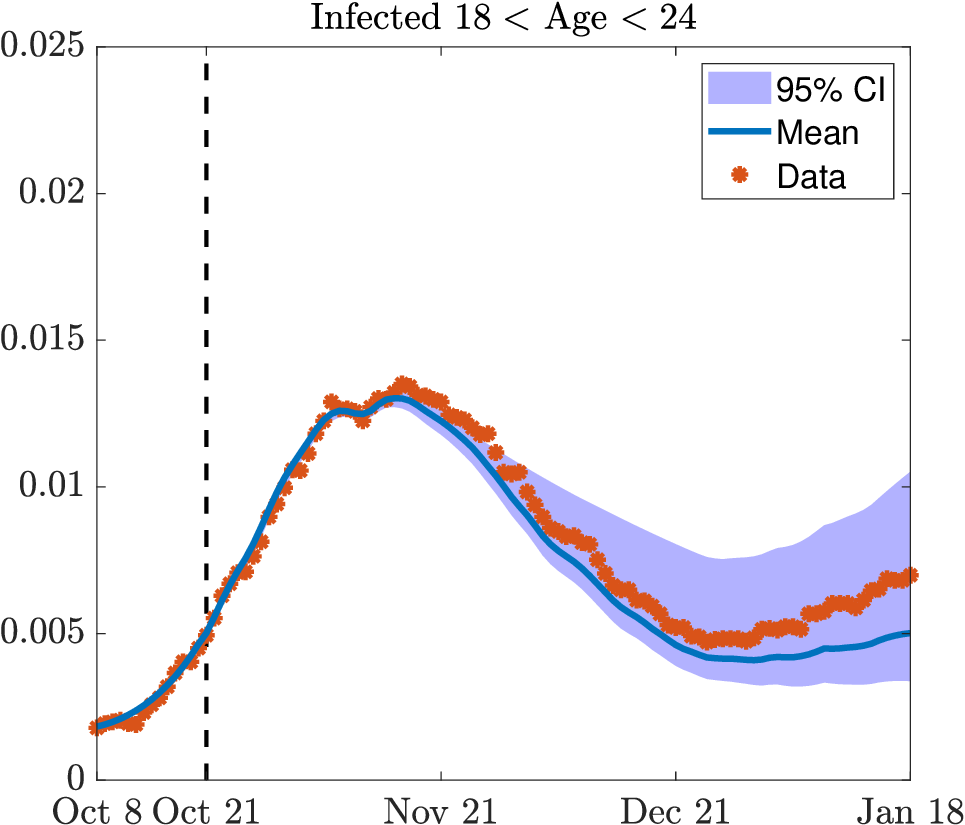}		\includegraphics[width=0.328\linewidth]{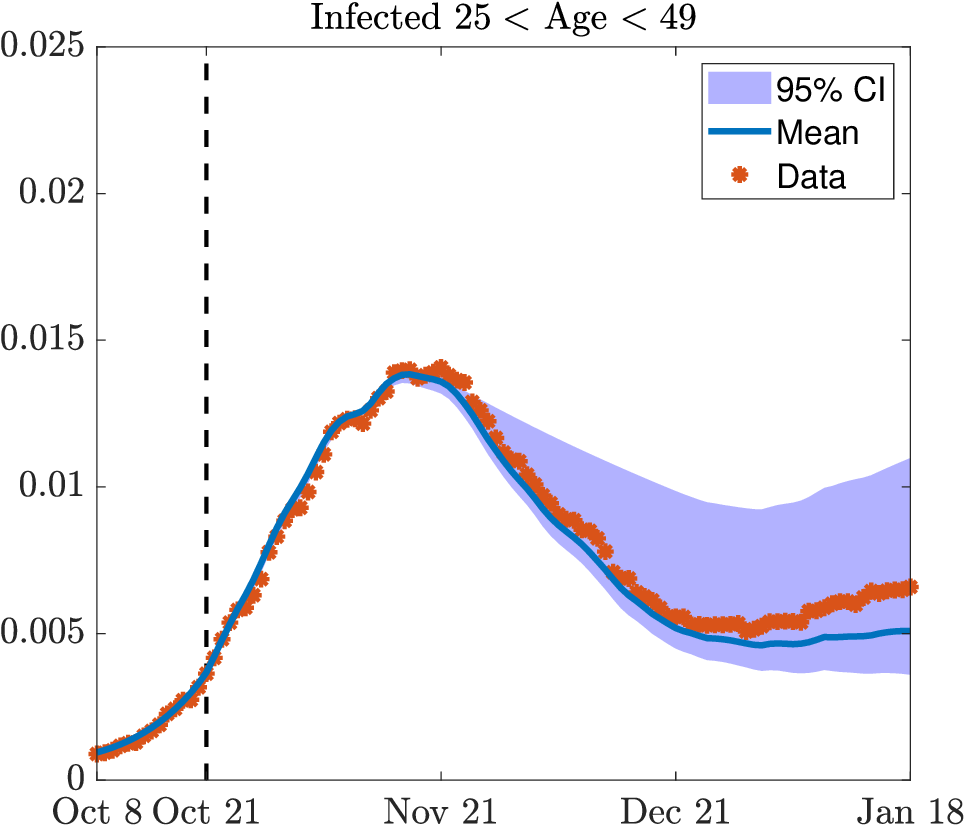}\\ \vspace{0.2cm}
	\includegraphics[width=0.328\linewidth]{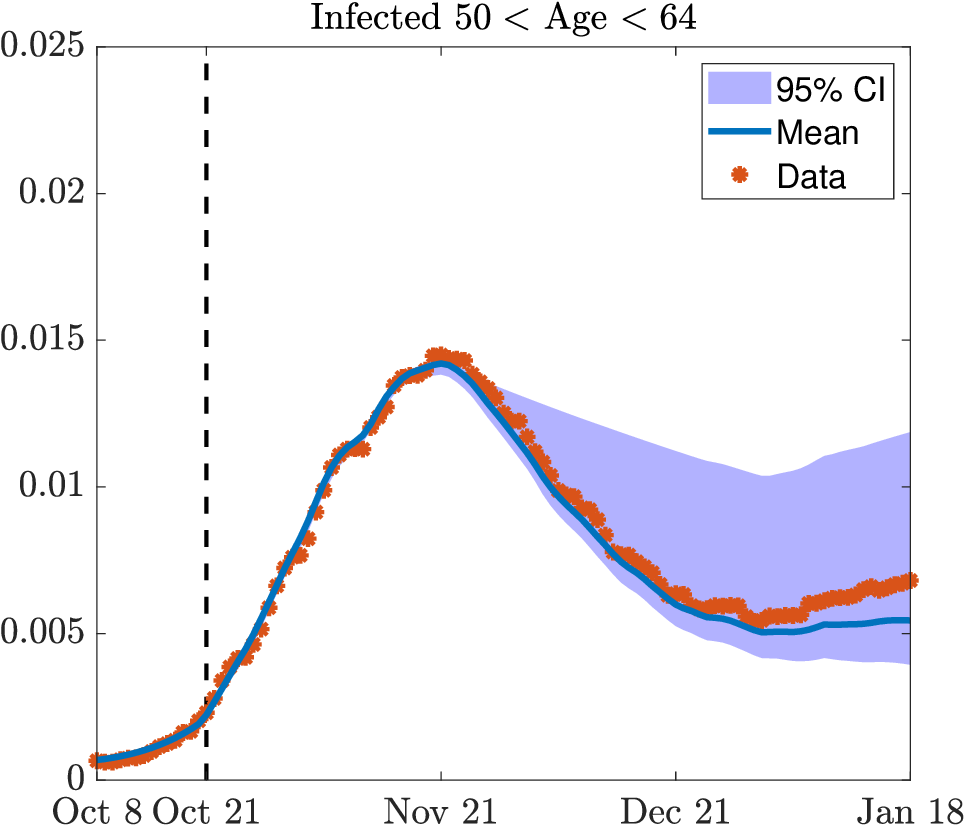}
	\includegraphics[width=0.328\linewidth]{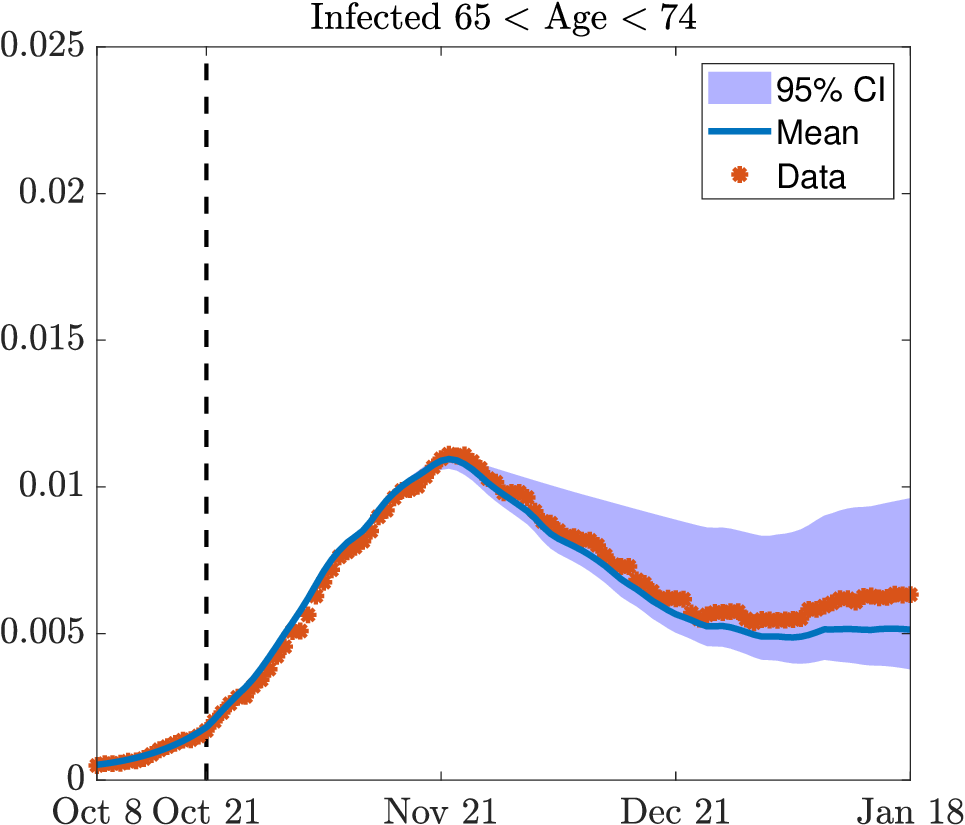}		\includegraphics[width=0.328\linewidth]{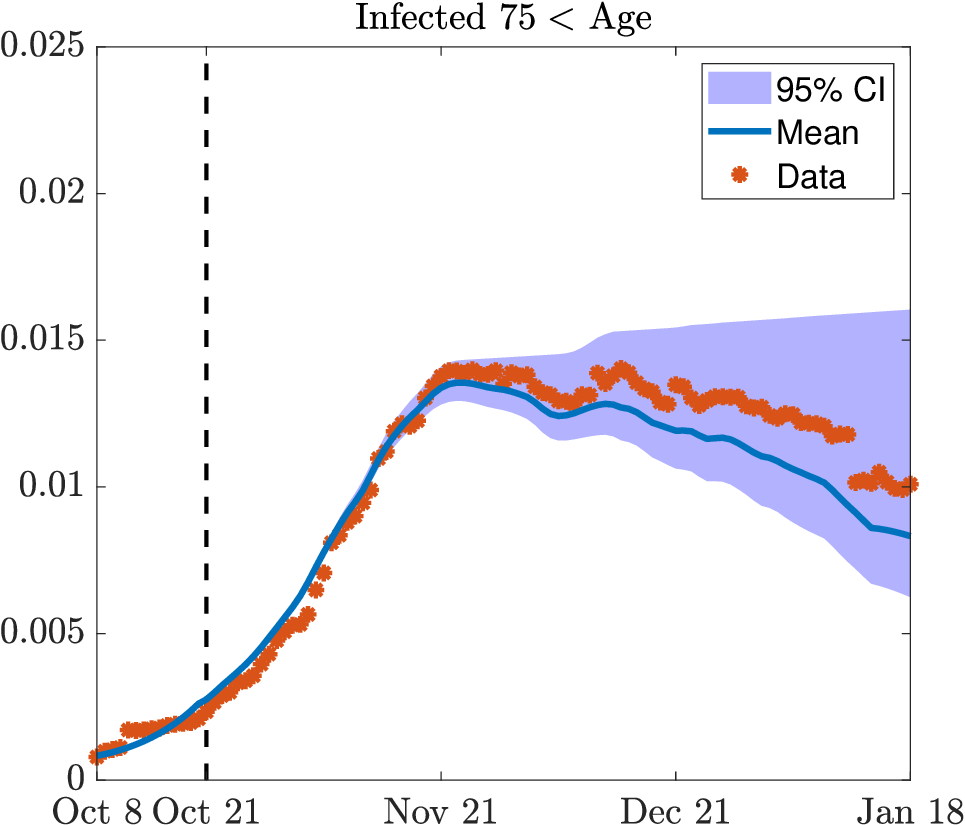}
	\caption{  Dynamics of the infected population obtained by solving the social-SIAR model \eqref{eq:social_SIAR_ages} and compared with experimental data. The plots show the mean epidemic trajectory with the 95\% confidence interval (shaded area), alongside the observed data. The black dashed line separates the two epidemic phases. Each image corresponds to a different age group. }
	\label{fig:infected_Ages}
\end{figure}

Once the epidemic parameters have been estimated, the second optimization procedure consists in identifying the optimal functions $H(x,t_j,z)$, which mimic the implementation of governmental interventions such as mask mandates and mobility restrictions, and $\xi(x,t_j,z)$, the number of unknown asymptomatic, for $t_j\in [t_j-k_l,t_j+k_r]$, with $j=1,\ldots,N_t$, where $N_t$ is the total number of considered time steps, $k_l=3$, and $k_r=4$. This choice corresponds to averaging the epidemic data over a one-week window, for each value of $z$ representing uncertainty, and for each age class $x$. The data used for this second optimization correspond to the subsequent phase of the second wave of the pandemic, specifically the interval $[t_L,T]$, with $t_L$ marking the start date (October 21st, 2020) and $T$ the end date (January 18th, 2021). As in the previous optimization, we sample the recovery rates from beta distributions and solve the following optimization problem for each $m=1,\ldots,M$ and for $t\in[t_j-k_l,t_j+k_r]$:
\begin{equation}\label{eq:min_H}
	\min_{H_{m}(x,t_j),\xi_{m}(x,t_j)} \int_{\mathcal{A}} \mathcal{J}(I(y,t,z_m),\hat{I}(y,t),R(y,t,z_m),\hat{R}(y,t)) d y,   
\end{equation}
where $\mathcal{J}(\cdot)$ is defined as in \eqref{eq:functional}. Both optimization problems \eqref{eq:min_beta}-\eqref{eq:min_H} have been solved by using the Matlab function \texttt{fmincon} combined with a RK4 integration method of the systems of ODEs. Figure \ref{fig:infected_noAges} shows the comparison of the above described optimization strategy for the social SIAR model \eqref{eq:social_SIAR} in terms of the number of infected individuals. In the figure, the mean trajectory and the 95\% confidence interval, together with the available reported data on infected cases are shown. The black dashed line distinguishes between first and second phase of the pandemic corresponding to the two optimization procedures. 
Figure \ref{fig:infected_Ages} shows the corresponding results for the age-structured model \eqref{eq:social_SIAR_ages}. As before, the mean trajectory and the 95\% confidence interval are displayed, together with the reported data on infected cases. The black dashed line separates the first and second phases of the pandemic. Each image corresponds to a different age group.  

\section{Predictions of COVID-19 dynamics using neural networks}\label{sec:neural_networks} 
Building upon the methodology introduced in our previous work \cite{awais2025data}, we aim to train neural networks on augmented datasets to enhance the accuracy of their solution both in terms of interpolations and predictions. In our earlier study, we considered deterministic synthetic data arising from the solution of a simplified epidemic model, namely the social-SIR model \eqref{eq:deterministicSIR_H}. Since the observational data were recorded on a daily basis (with time step $h = 1$), we generated additional data over the same interval by solving the system of ODEs with a finer resolution ($h = 0.2$).  In this work, we extend this strategy to more complex models \eqref{eq:social_SIAR}-\eqref{eq:social_SIAR_ages}, with the important distinction that we now incorporate uncertainty in the model parameters, as outlined in the previous sections. Our goal is to train neural networks that are well-suited for predictive tasks, such as Physics Informed Neural Networks (PINNs), and Non-Linear Autoregressive (NAR) networks. Both PINNs and NAR networks are implemented using a Feed-Forward neural network architecture with $L+1$ layers, as follows
\begin{equation}\label{eq:NN_structure} 
	\begin{split}
		&x^1 = W^1 x + b^1, \qquad x^{l} = \sigma \circ (W^lx^{l-1} + b^l), \quad l=2,\ldots,L-1,\\
		&f^{NN}(x;{\bf b},{\bf W} ) = W^L z^{L-1}+ b^L,
	\end{split}
\end{equation} 
where $x = (x^1,\ldots,x^d) \in \mathbb{R}^d$ is the input, $\sigma$ is an activation function  to convert the input signals to output signals, and ${\bf W} = (W^1,\ldots, W^L)$ and ${\bf b} = (b^1,\ldots,b^L)$ represent the weights and biases.  The weights and biases associated with these connections serve as the parameters of the network, adjusted iteratively through techniques like gradient descent or Adam method during the training phase to minimize prediction errors and enhance model performance. To find the optimal parameters $ \theta^* = (W^1,b^1,\ldots,W^L,b^L)$, we solve the following minimization problem
\begin{equation}\label{eq:min_prob}
	\theta^* = \arg \min_{\theta} \mathcal{L}(\theta),
\end{equation}
being 
\begin{equation}\label{eq:loss}
	\mathcal{L}(\theta) = \omega_d \mathcal{L}_d(\theta) + \omega_p \mathcal{L}_p (\theta),
\end{equation}
the loss function, for $\omega_d\geq 0$, $\omega_p\geq 0$. Specifically, in \eqref{eq:loss}, the function $\mathcal{L}_d(\theta)$ measures the discrepancy between the neural network solution and the data, while $ \mathcal{L}_p(\theta)$ encodes the physics of the problem. 
In the following, we present PINNs and NAR networks in the context of the age-structured model \eqref{eq:social_SIAR_ages}. A similar formulation applies to the simpler model \eqref{eq:social_SIAR}, where integration is performed over all age classes  $x\in \mathcal{A}$, with $\mathcal{A}$ denoting the age-classes.

\subsection{Physics Informed Neural Networks}\label{sec:PINNs}
PINNs offer a physics-driven learning framework where the network is trained not only to fit the data but also to satisfy the underlying differential equations governing the system dynamics. This makes them particularly appealing in scenarios where data are scarce or noisy, as the physical constraints serve as an effective regularizer. In our context, we define the PINN using the epidemic model equations, incorporating uncertain parameters $z$ sampled as described earlier. The architecture consists of:
\begin{itemize}
	\item \textbf{Input Layer:} Takes as input both the age-time variables $x,t$ and the uncertainty parameter $z$.
	\item \textbf{Hidden Layers:} Several fully connected layers with non-linear activation functions (e.g., tanh) to approximate the solution. 
	\item \textbf{Output Layer:} Returns the approximated solution values  $f_i^{NN}(x,t,z)$ for $i\in \mathcal{I} =\{S,I,A,R\}$ for any $x$, $t$ and $z$.
\end{itemize} 
To determine the optimal parameters, we solve the minimization problem in \eqref{eq:min_prob}, defining the loss function as in \eqref{eq:loss} with $\omega_d,\omega_p>0$. Specifically, the loss function encoding the data constraints reads as follows 
\begin{equation}
	\mathcal{L}_d(\theta) =  \sum_{i\in \mathcal{I}} \sum_{n=1}^{N_c}  \int_{\mathcal{A}} \left( \bar{f}_i^{NN}(x,t_n;\theta) - \bar{f}_i(x,t_n) \right)^2 dx,
\end{equation}
being 
\begin{equation}\label{eq:mean_sol}
	\bar{f}^{NN}_i(x,t_n;\theta) = \sum_{m=1}^M f^{NN}_i(x,t_n,z_m;\theta) w_m,\qquad 	\bar{f}_i(x,t_n;\theta) = \sum_{m=1}^M \hat{f}_i(x,t_n,z_m;\theta) w_m, 
\end{equation}
where $w_m$ are the weights associated to the uncertainty values $z_m$, and $N_c$, $M$ are the total number of samples that we consider in time and in the uncertainty space respectively, while $\hat{f}_i$ are the data representing the number of Susceptible, Infected, Asymptomatic, and Recovered individuals.

Via automatic differentiation we then compute the derivative with respect to time of the quantities $f^{NN}_i(x,t,z;\theta)$ for any $i\in \mathcal{I}$, and we define the physical loss as 
\begin{equation}\label{eq:loss_physics}
	\mathcal{L}_p(\theta) = \sum_{i\in \mathcal{I}} \sum_{n=1}^{N_c}  \int_{\mathcal{A}}  \mathcal{\bar{R}}^2_i(x,t_n;\theta)    dx,
\end{equation} 
being 
\begin{equation}
	\mathcal{\bar{R}}_i(x,t_n;\theta) =  \sum_{m=1}^M \mathcal{R}_i(x,t_n,z_m;\theta) w_m, 
\end{equation}
with  $w_m$ weights associated with the Gauss-Jacobi nodes, and 
\begin{equation}\label{eq:residuals}
	\begin{split}
		&\mathcal{R}_S  = \partial_t f^{NN}_S(x,t_n,z_m;\theta) + \Lambda(x,t_n,z_m;\theta),\\
		&\mathcal{R}_I  = \partial_t f^{NN}_I(x,t_n,z_m;\theta) - \xi(x,t_n,z_m) \Lambda(x,t_n,z_m;\theta) + \gamma_{I}(x,z_m) f^{NN}_I(x,t_n,z_m;\theta), \\
		&\mathcal{R}_A = \partial_t f^{NN}_A(x,t_n,z_m;\theta) - (1-\xi(x,t_n,z_m)) \Lambda(x,t_n,z_m;\theta) + \gamma_{A}(x,z_m) f^{NN}_A(x,t_n,z_m;\theta),\\
		&\mathcal{R}_R  = \partial_t f^{NN}_R(x,t_n,z_m;\theta)-  \gamma_{I}(x,z_m) f^{NN}_I(x,t_n,z_m;\theta) -  \gamma_{A}(x,z_m) f^{NN}_A(x,t_n,z_m;\theta),
	\end{split}
\end{equation}
where 
\begin{equation}\label{eq:Lambda_NN}
	\begin{split} 
		\Lambda  = &\beta(x,z_m) f_S^{NN}(x,t_n,z_m;\theta) H(x,t_n,z_m)  \\ & \qquad \qquad \qquad \int_{\mathcal{A}} H(y,t_n,z_m) (kf_I^{NN}(y,t_n,z_m;\theta) + f_A^{NN}(y,t_n,z_m;\theta) ) dy ,
	\end{split}
\end{equation}
for $k\in[0,1]$. 
\subsection{Non-Linear Autoregressive networks}\label{sec:NARs}
NAR networks leverage past time series values to forecast future states, making them particularly effective for time-dependent predictions, especially in short term forecasting. The NAR architecture includes the following components:
\begin{itemize}
	\item \textbf{Input Layer:} Composed of infected population values at previous time steps, $I(x,t-d,z),\ldots,I(x,t-1,z)$, for a chosen delay $d\geq1$ and for any value of the uncertain variable $z$ and of the age variable $x$.
	\item \textbf{Hidden Layers:} A set of fully connected layers equipped with nonlinear activation functions (e.g., ReLU or tanh) to model complex temporal dependencies.
	\item \textbf{Output Layer:} Produces the predicted value $I(x,t,z)$, representing the number of infected at time $t$ for each realization of $z$ and of the age variable $x$.
\end{itemize}
Once trained, the NAR network can recursively predict the epidemic evolution by feeding its own predictions as inputs for future time steps (closed-loop strategy). The NAR approach is particularly effective for qualitative short-term forecasting, as it can directly learn epidemic dynamics from data without relying on the underlying model equations. This data-driven nature allows NAR to outperform PINNs in the short-term setting, both in terms of accuracy and computational efficiency. However, when it comes to long-term forecasting, PINNs provide a more reliable framework, as their physics-informed structure enables them to incorporate the underlying epidemic dynamics and capture long-term trends that NAR can not.\\ 
The training process involves the minimization of the loss function \eqref{eq:loss} where we suppose $\omega_d>0$, $\omega_p=0$, and 
\begin{equation}\label{eq:loss_NAR} 
	\mathcal{L}_d(\theta) =  \sum_{n=1}^{N_c}\sum_{m=1}^{M}\int_{\mathcal{A}}  \left( I^{NN}(x,t_n,z_m;\theta) - \hat{I} (x,t_n,z_m) \right)^2 dx,
\end{equation}
being $I^{NN}$ the solution computed by the neural network and $\hat{I}$ the data, with $N_c$ denoting the number of samples in time, and $M$ the total number of samples in the uncertainty space.

\begin{remark}\label{rmk:remark1} 
	In the case of PINNs, it is essential to include constraints for all population compartments in the loss function to accurately capture the dynamics prescribed by the social SIR model, even if the primary objective is to reproduce only the infected population. In contrast, for the NAR network, the dynamics are not explicitly enforced in the loss function but instead the calibrated model is used to produce an augmented data set used to feed the network. 
\end{remark}

\section{Numerical experiments}\label{sec:numerical_exp} 
We now proceed with numerical experiments to validate our methodology. In particular, we aim to show that neural networks trained on larger synthetic datasets achieve higher accuracy in both prediction and extrapolation tasks. In addition, the goal is to demonstrate that NAR networks provide a viable alternative to PINNs in the context of short term predictions, especially when trained on synthetic data.
The synthetic datasets are generated by solving systems \eqref{eq:social_SIAR}–\eqref{eq:social_SIAR_ages} over the time interval $[t_0, T]$, with $t_0 = 15$ and $T = 105$, using a time step of $h = 0.2$. This time window corresponds to the second phase of the pandemic in Italy, specifically from October 21st, 2020 to January 18th, 2021.  We empathize that data are generated over the entire time horizon to guarantee sufficient samples for both short- and long-term forecasting. However, the neural networks are trained exclusively on data within the specified training window, without any use of information from the test period. Model parameters are treated as uncertain, as previously described, and uncertainty is represented by sampling $M = 5$ Gauss-Jacobi nodes $z_m$.  The macroscopic incidence rates $H_J(\cdot)$ for $J\in\{S,I,A\}$, and the number of asymptomatic $\xi(\cdot)$ are reconstructed from the discrete dataset, computed by solving the minimization problem \eqref{eq:min_H}, by selecting the value corresponding to the time index immediately preceding the desired time $t$.  For sake of comparison, we also train the same neural networks only using the real data at disposal. 
We will mainly focus on the short term forecast, splitting the datasets into training and test sets: the training set covers the period from October 21st to January 8th, while the test set spans January 9th to January 18th. At the end of the section, we compare the performances of PINNs and NARs networks focusing on the short and long term forecasting. For the long term forecasting we split the dataset into training set (from October 21st to November 19th) and test set (from November 20th to January 3rd), in order to capture the peak of the pandemic.

\subsection{Physics informed neural network.}\label{sec:PINNs_numerics}
We start by defining two Feed-Forward network architectures  with width 32 and 3 hidden layers, using \texttt{tanh} as activation function. The networks are designed to approximate the solutions of the social SIAR model and the age-structured social SIAR model by means of PINNs. Both networks are trained on real and synthetic data, assuming in the loss function \eqref{eq:loss}, $\omega_d=\omega_p=1$ to account both for the data and the underline physics. For real data, in the case of the non-age structured model, the trained network takes the time $t$ as input and it produces as  output the corresponding value $f_i(t)$ for any compartment $i\in\{S,I,A,R\}$ at time $t$, for any $t$ corresponding to real data. In the case of synthetic data, the input consists of both the time $t$ and the uncertainty variable $z$, and the output is the value $f_i(t,z)$ for each compartment $i\in\{S,I,A,R\}$, for any $t$ corresponding to synthetic data. In the age-structured setting, the network input is the triplet $(x,t,z)$ and the output is $f_i(x,t,z)$, where $x\in \mathcal{A}$, with $\mathcal{A}$ denoting the different age classes.
The networks are trained for 50000 epochs using the Adam optimizer with a learning rate $10^{-2}$.  Following the procedure described previously, we ensure that both the data and the physics of the system are learned by the network. 

\begin{figure}[h] 
	\centering
	\includegraphics[width=0.426\linewidth]{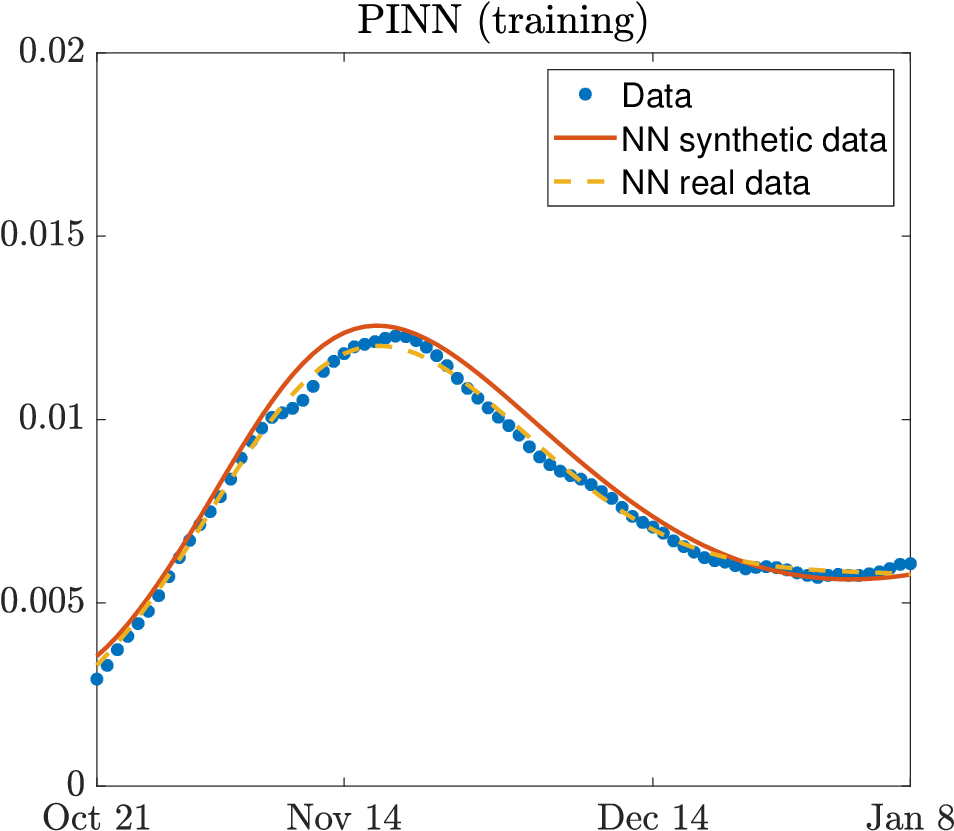}
	\includegraphics[width=0.415\linewidth]{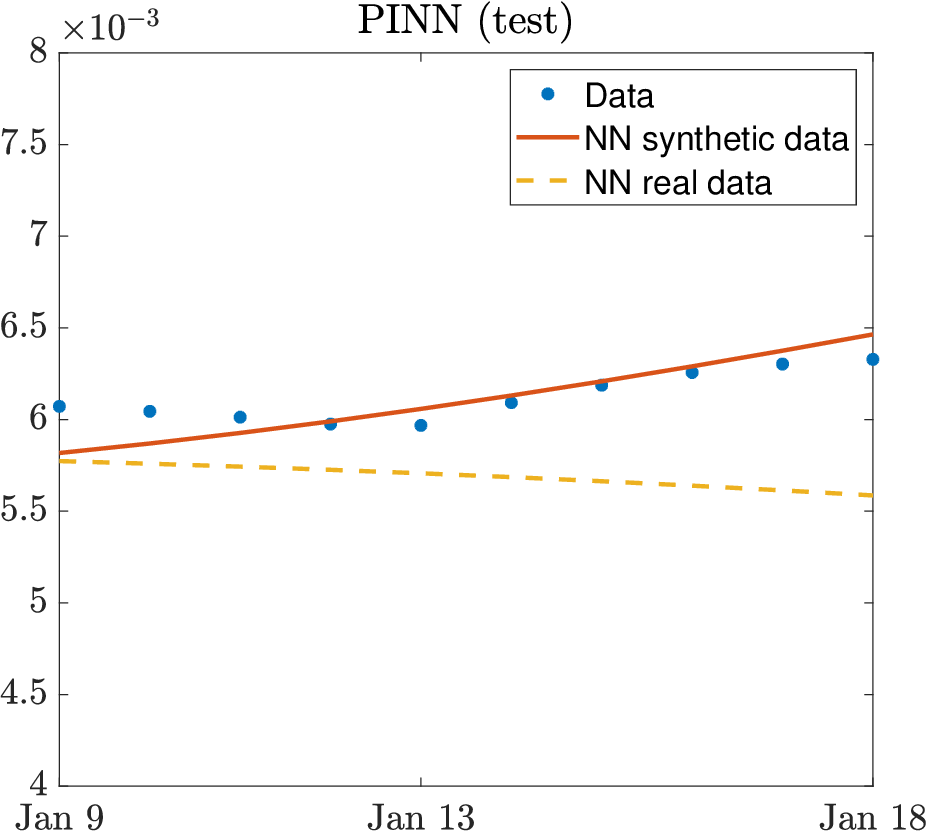} 
	\caption{Physics informed neural network for the social SIAR model \eqref{eq:social_SIAR}. Solution obtained by training a PINN network on both real and synthetic data, compared to the available data. On the left, the solution computed on the training set. On the right, the solution computed on the test set.}
	\label{fig:PINN_noAges}
\end{figure}

\begin{figure}[h] 
	\centering
	\includegraphics[width=0.328\linewidth]{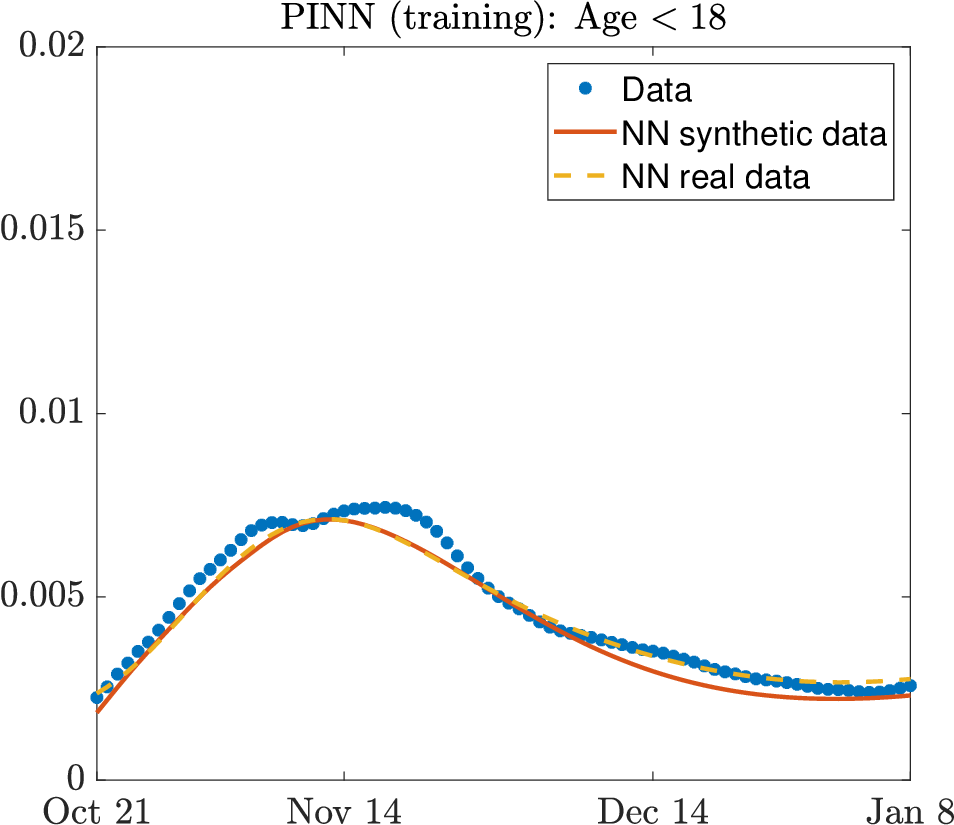}
	\includegraphics[width=0.328\linewidth]{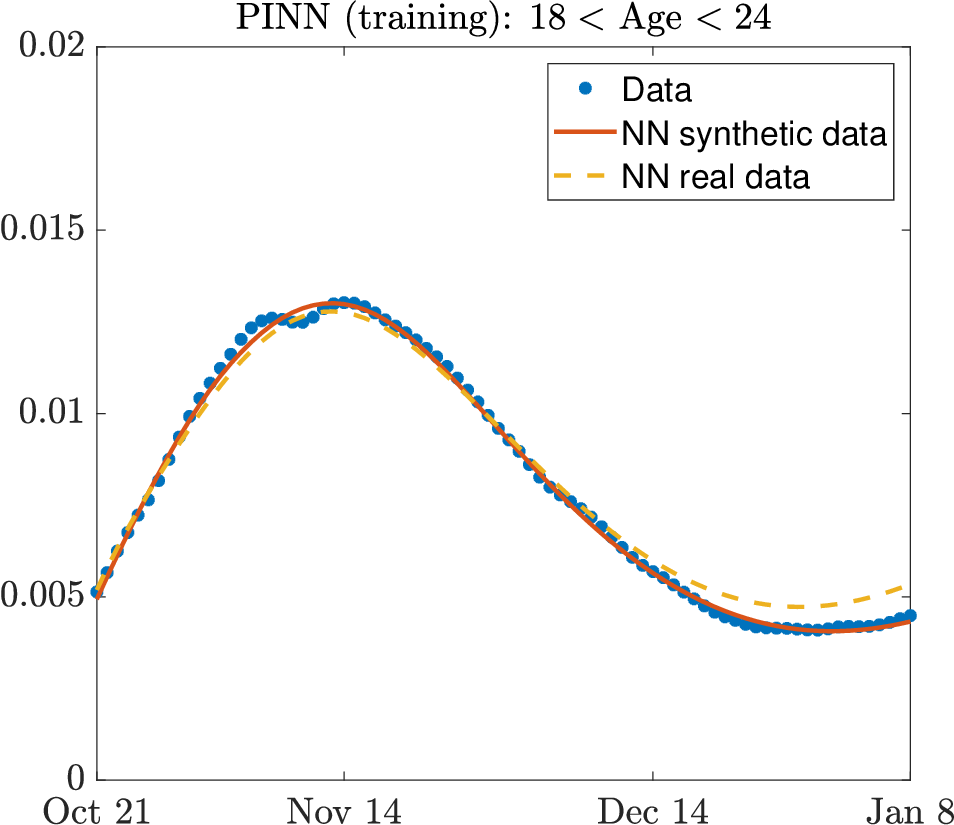}
	\includegraphics[width=0.328\linewidth]{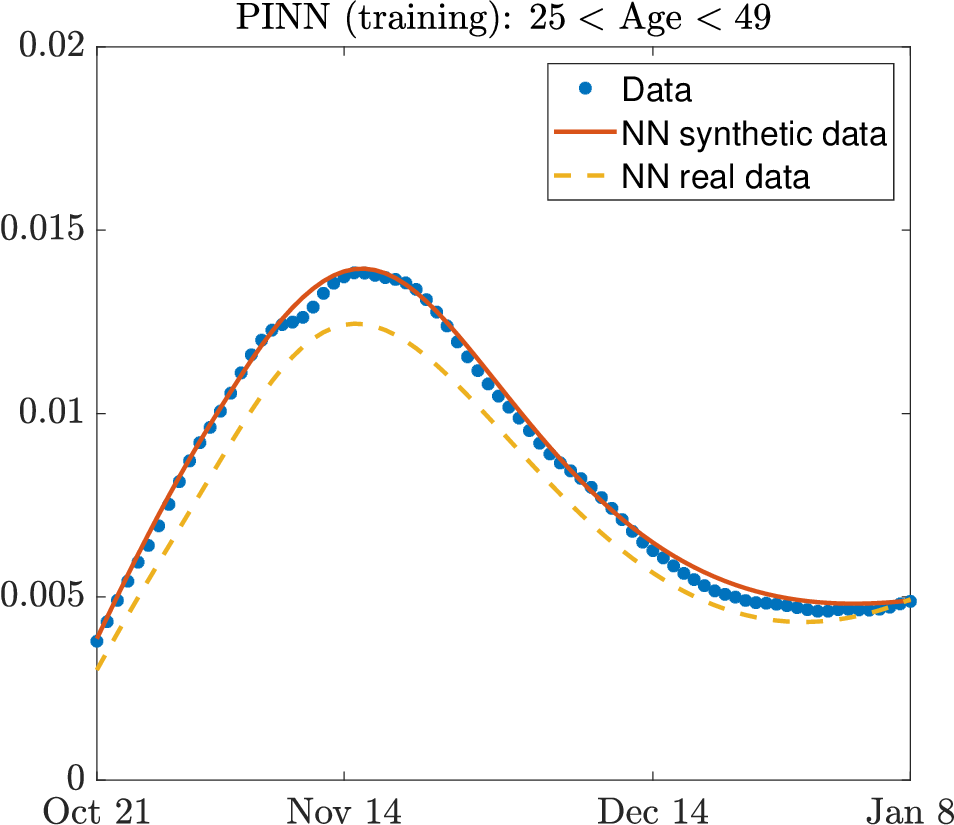}\\\vspace{0.2cm}
	\includegraphics[width=0.328\linewidth]{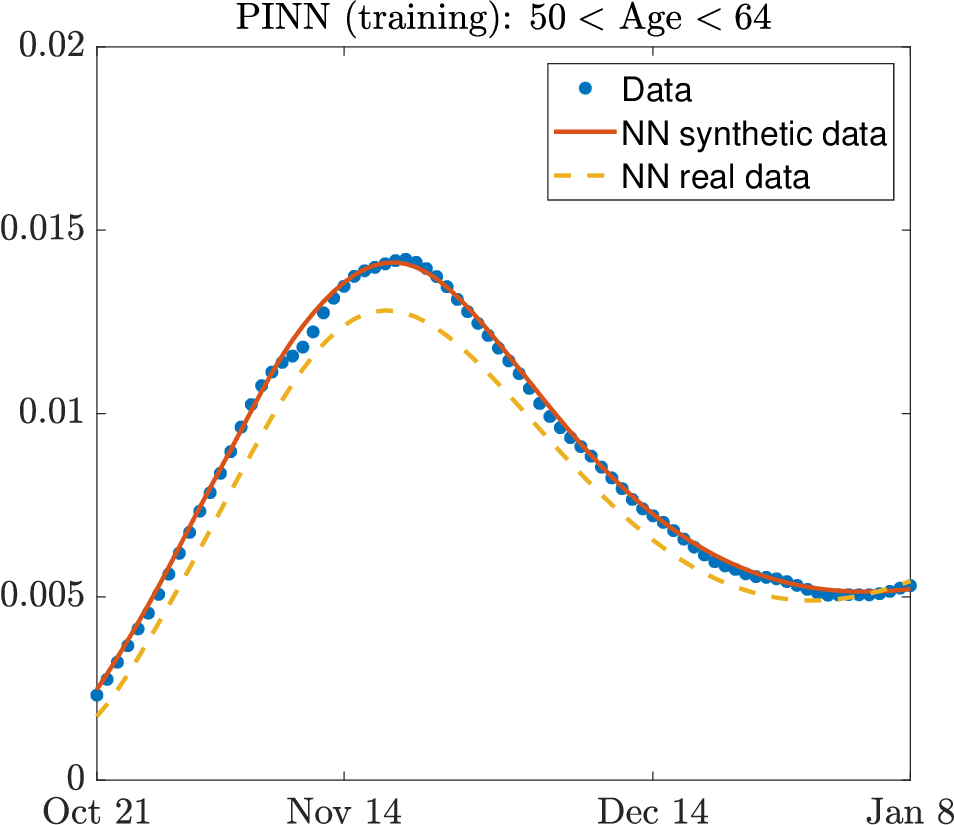}
	\includegraphics[width=0.328\linewidth]{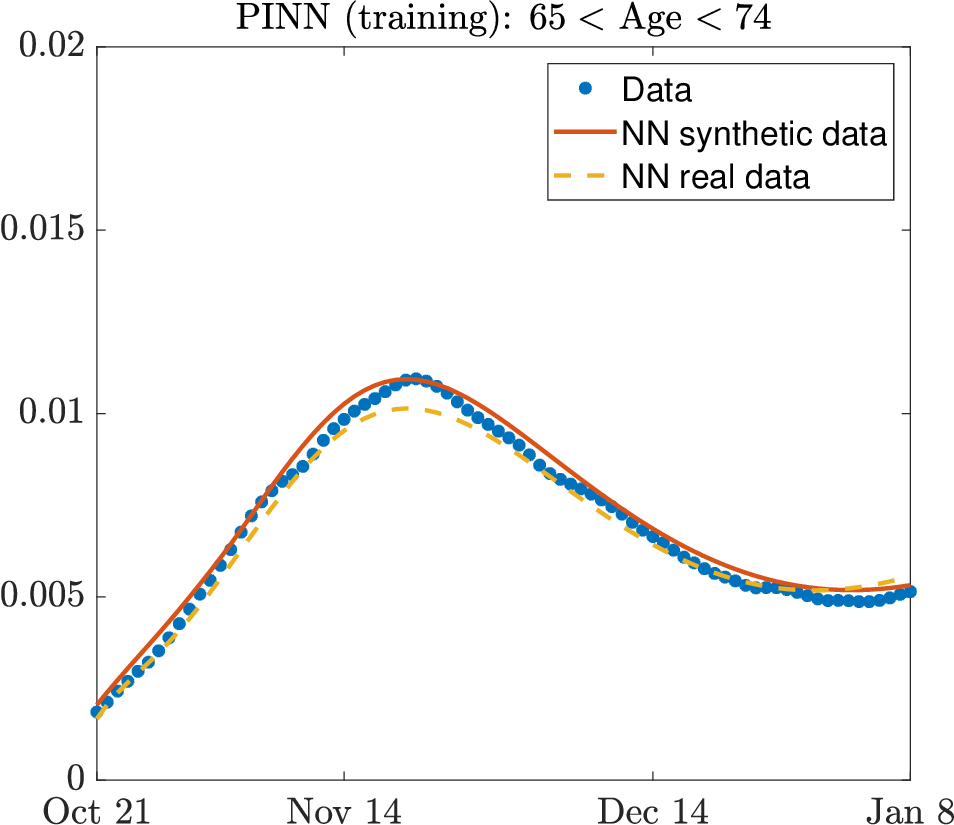}
	\includegraphics[width=0.328\linewidth]{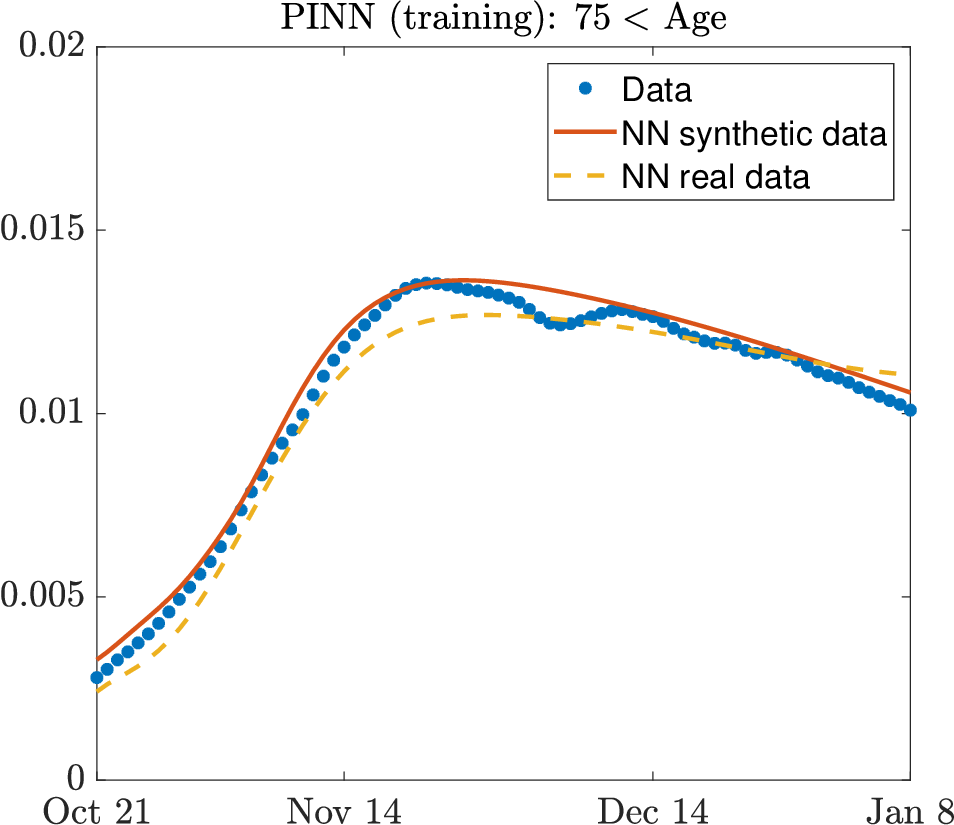}
	\caption{Physics informed neural network (training set) for the age-structured social SIAR model \eqref{eq:social_SIAR_ages}. Solution obtained by training a PINN network on both real and synthetic data, compared to the available data. Each plot corresponds to a different age class. Training set.}
	\label{fig:PINN_ages_training}
\end{figure}

\begin{figure}[h] 
	\centering
	\includegraphics[width=0.328\linewidth]{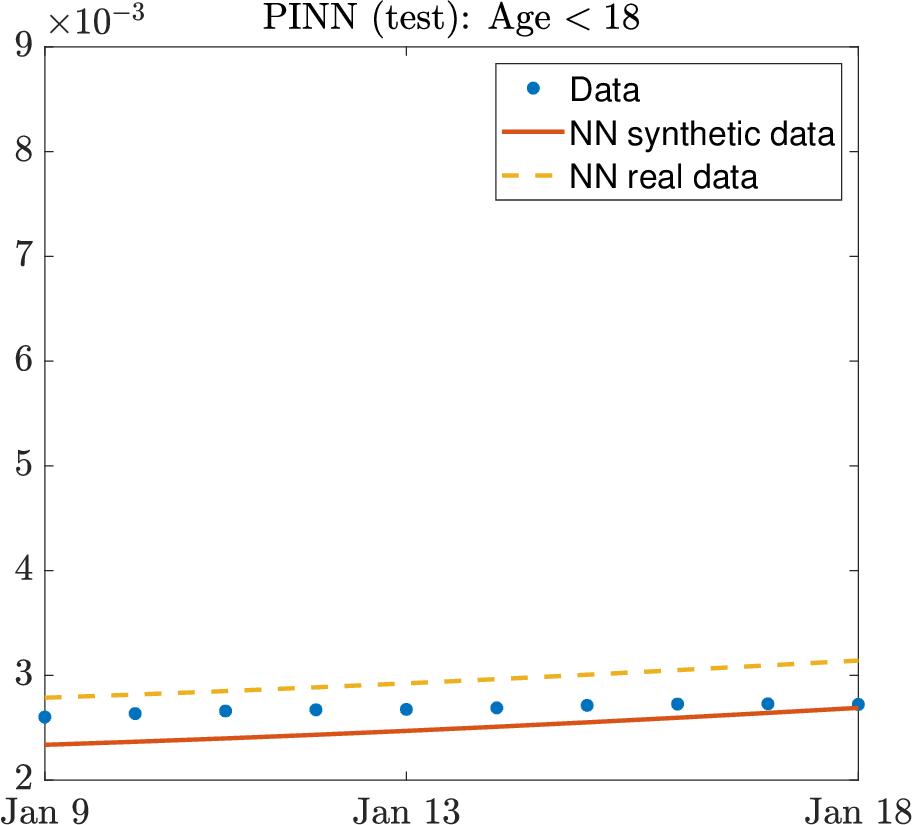}
	\includegraphics[width=0.328\linewidth]{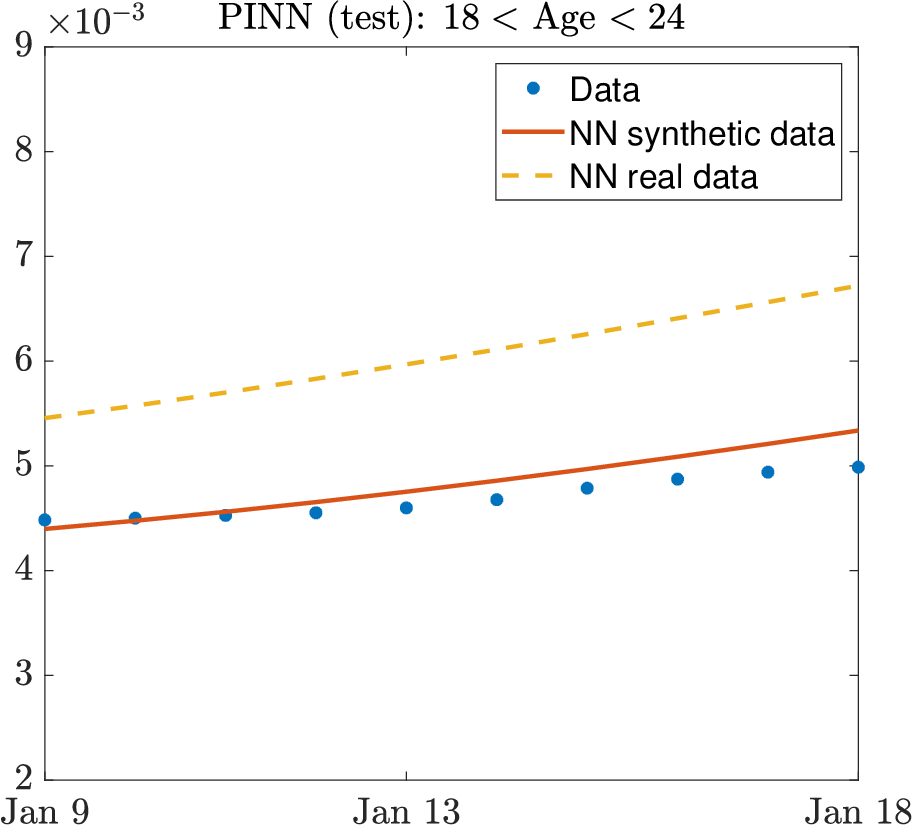}
	\includegraphics[width=0.328\linewidth]{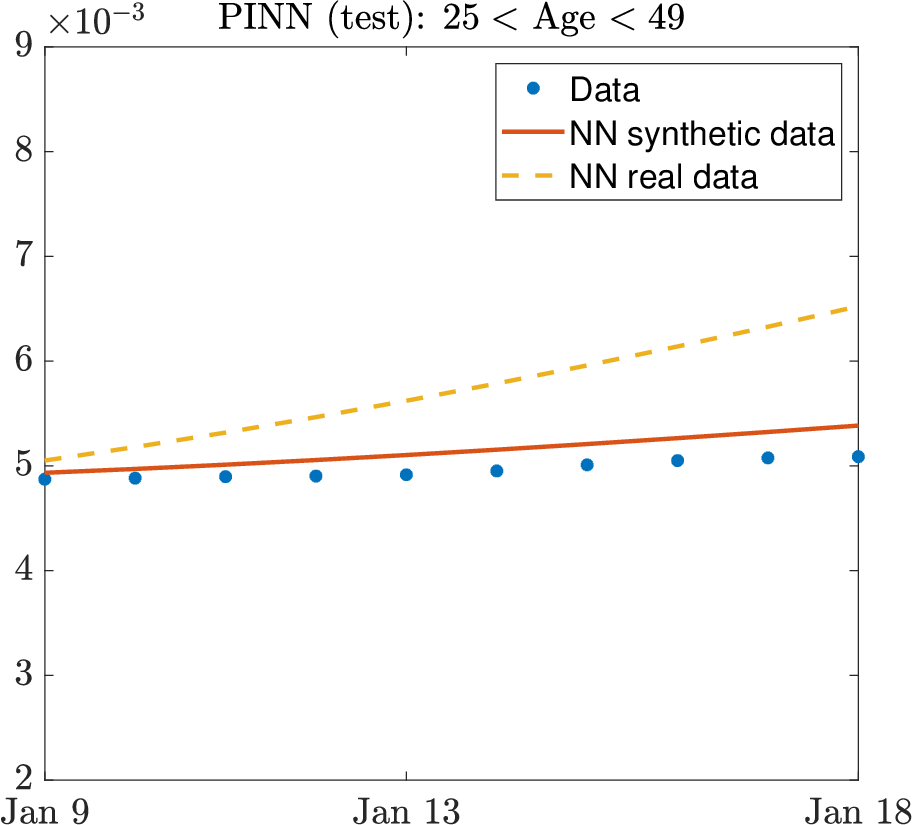}\\\vspace{0.2cm}
	\includegraphics[width=0.328\linewidth]{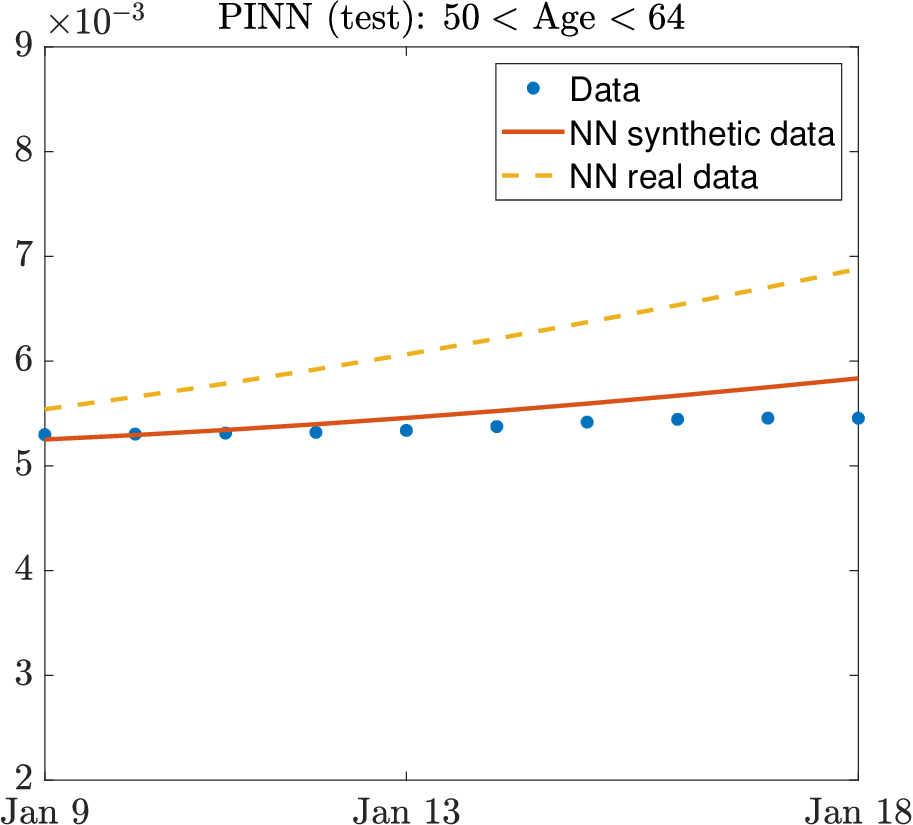}
	\includegraphics[width=0.328\linewidth]{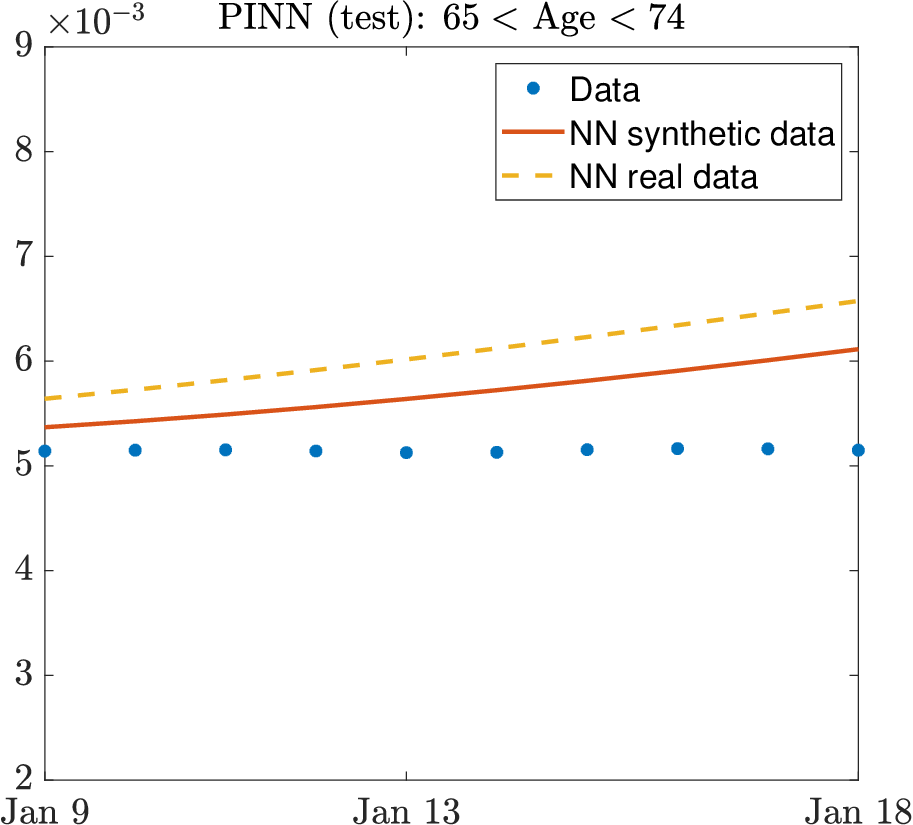}
	\includegraphics[width=0.329\linewidth]{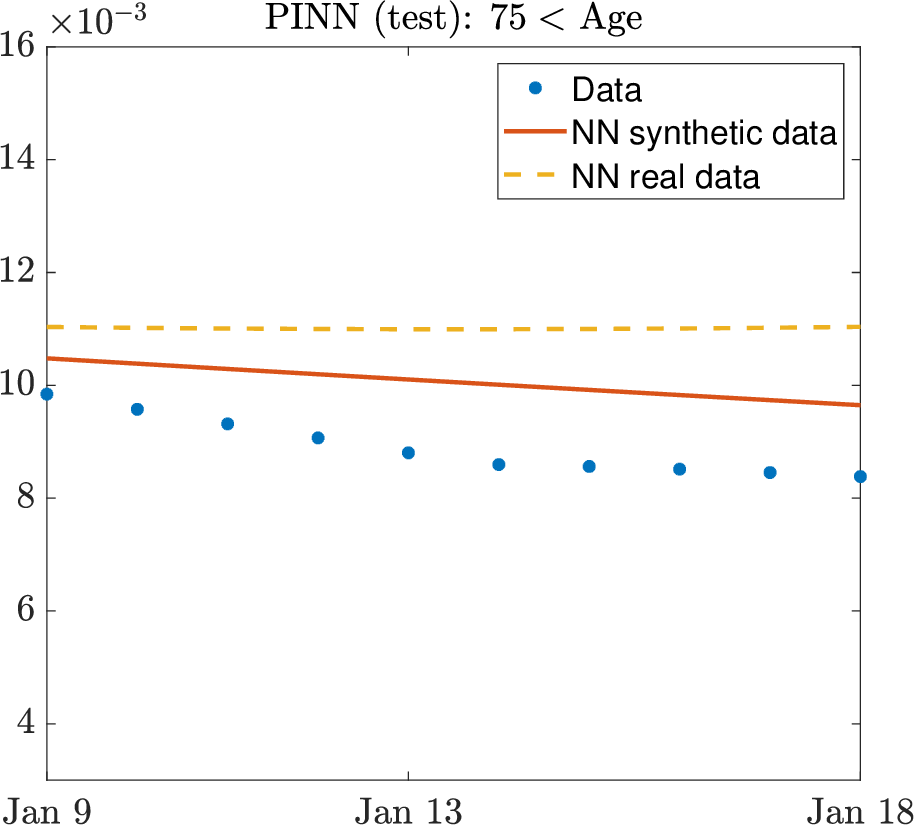}
	\caption{Physics informed neural network (test set) for the age-structured social SIAR model \eqref{eq:social_SIAR_ages}. Solution obtained by training a PINN network on both real and synthetic data, compared to the available data. Each plot corresponds to a different age class. Test set.}
	\label{fig:PINN_ages_test}
\end{figure}

In Figure \ref{fig:PINN_noAges} the comparison between the PINN solutions and the available data in the case of the social SIAR model \eqref{eq:social_SIAR} is reported. On the left, the image shows the part relative to the training set while on the right the image refers to the test set. We clearly see that the network trained on real data achieves a better fit on the training set, whereas the network trained on synthetic data produces a more accurate solution on the test set and thus it is able to better forecast the time evolution of the epidemic. 
Figure \ref{fig:PINN_ages_training} shows the comparison between the PINN solutions and the available training data in the case of the age-structured social SIAR model \eqref{eq:social_SIAR_ages} for the training part. Figure \ref{fig:PINN_ages_test} shows the same results but computed on the test set. In most of the cases, despite of the complexity of the solution, the networks achieve a good approximation in interpolation; however, they still face challenges in prediction accuracy, even when synthetic data are used to enhance their forecasting capabilities.
\subsection{Non-linear autoregressive network.} \label{sec:NARs_numerics}
\begin{figure}[h] 
	\centering
	\includegraphics[width=0.426\linewidth]{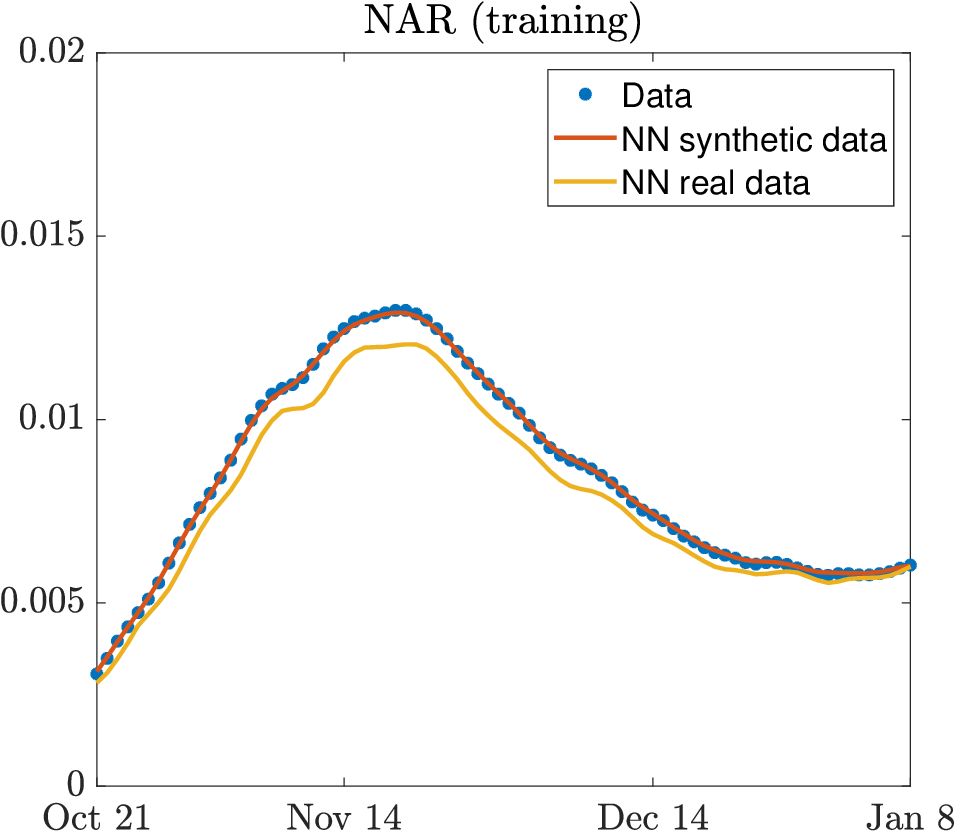}
	\includegraphics[width=0.415\linewidth]{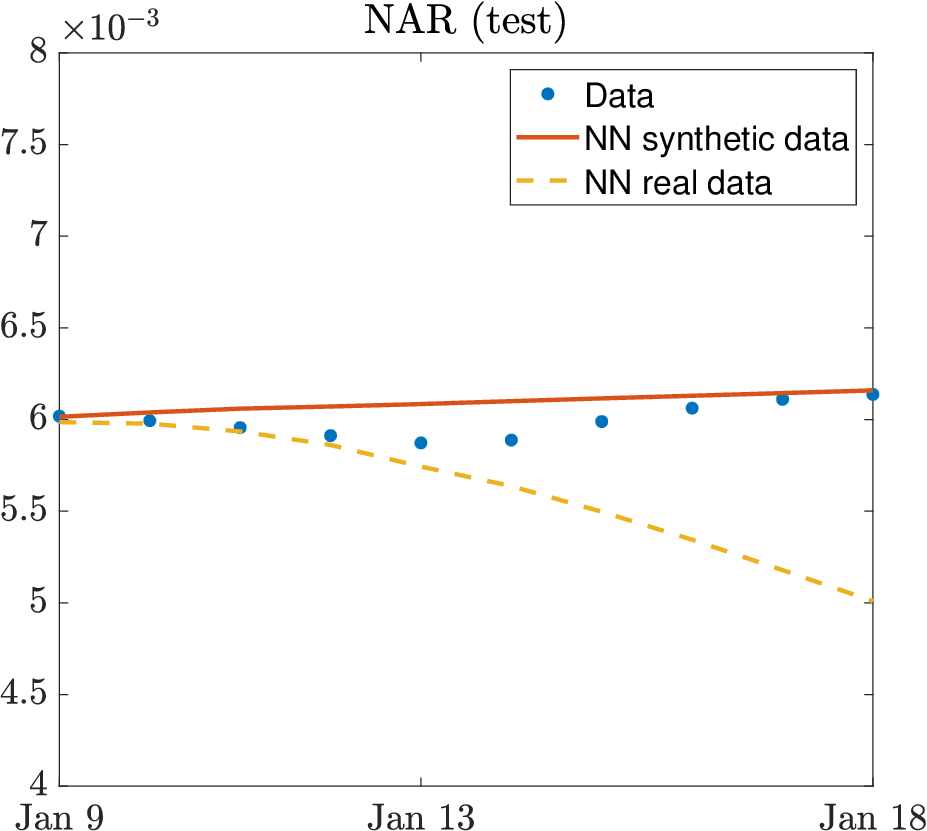} 
	\caption{Nonlinear autoregressive network for the social SIAR model \eqref{eq:social_SIAR}. Solution obtained by training a NAR network on both real and synthetic data, compared to the available data. On the left, the solution computed on the training set. On the right, the solution computed on the test set.}
	\label{fig:NAR_noAges}
\end{figure}

\begin{figure}[ht] 
	\centering
	\includegraphics[width=0.328\linewidth]{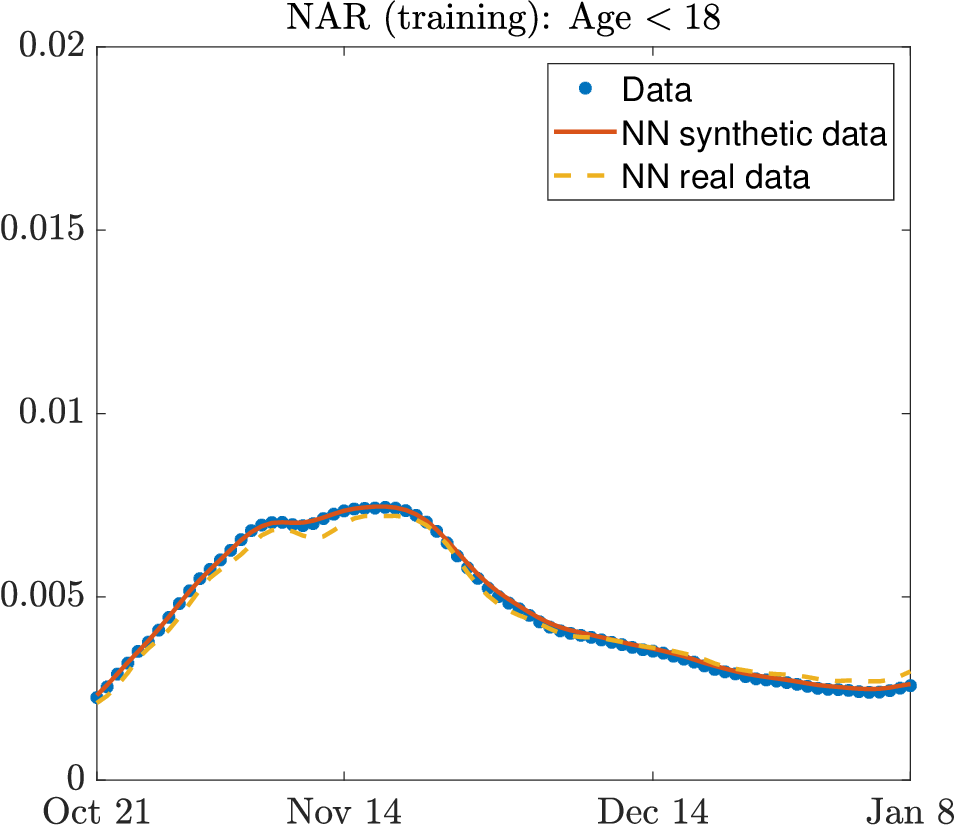}
	\includegraphics[width=0.328\linewidth]{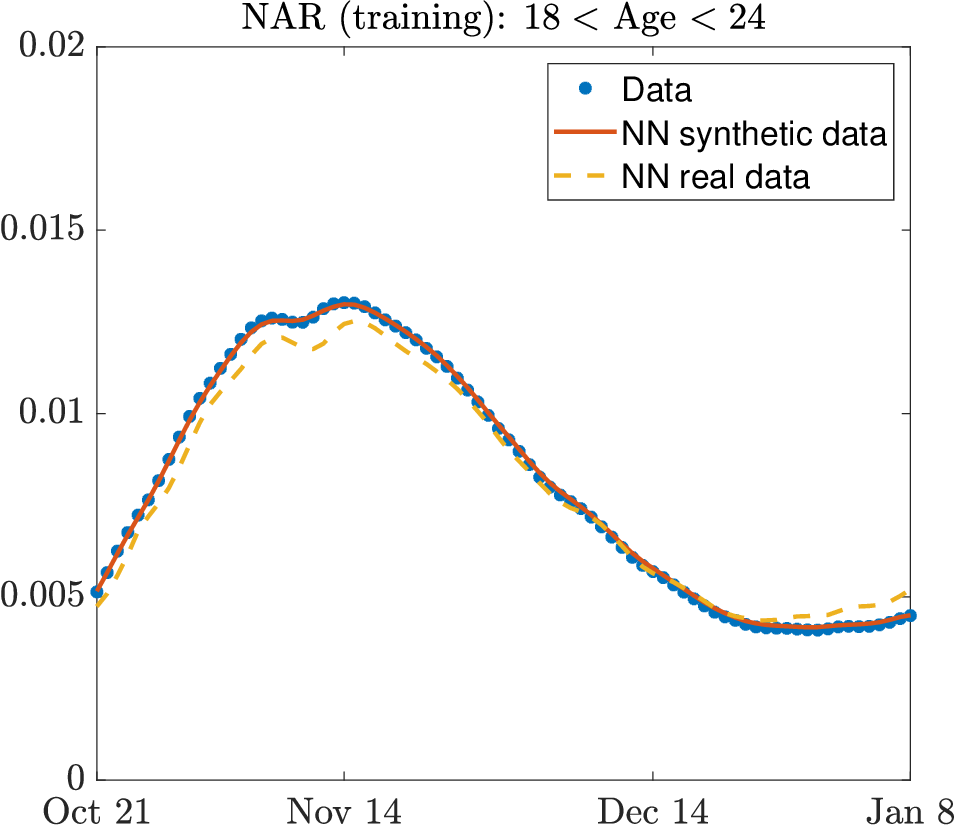}
	\includegraphics[width=0.328\linewidth]{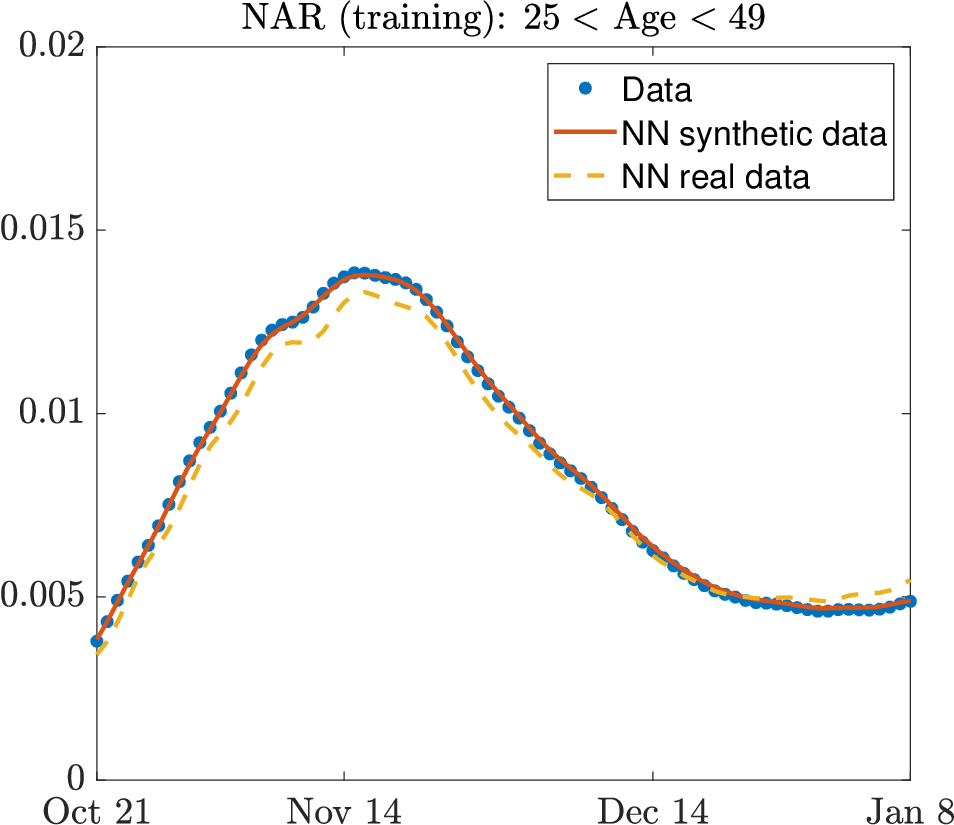}\\\vspace{0.2cm}
	\includegraphics[width=0.328\linewidth]{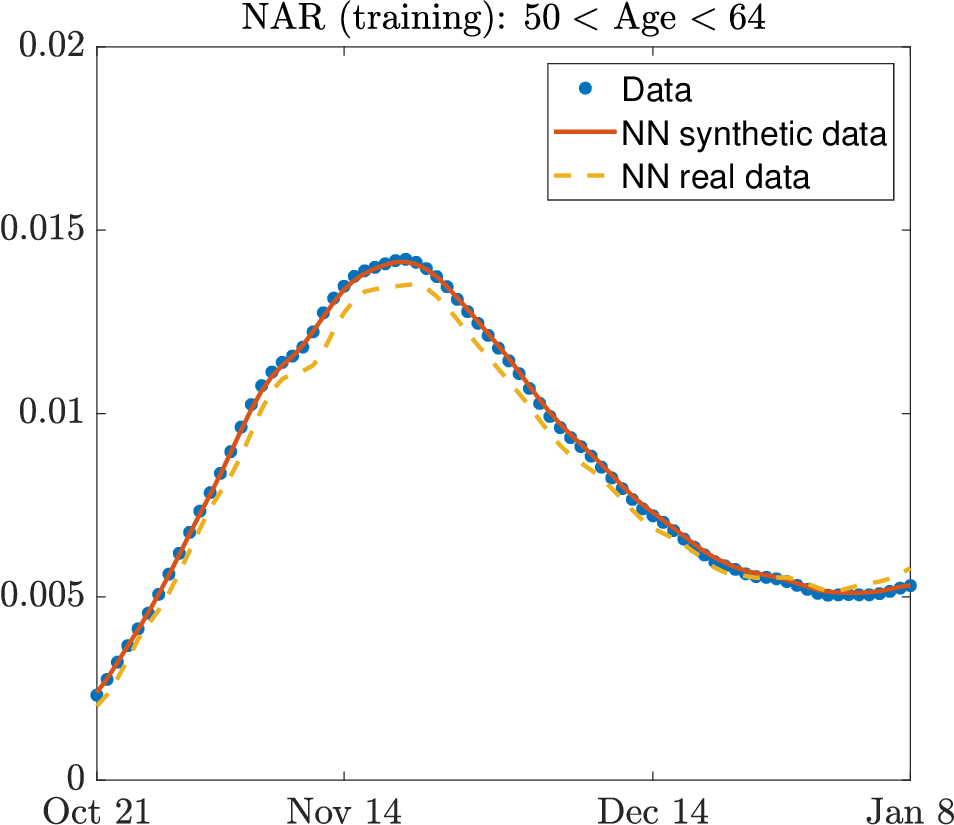}
	\includegraphics[width=0.328\linewidth]{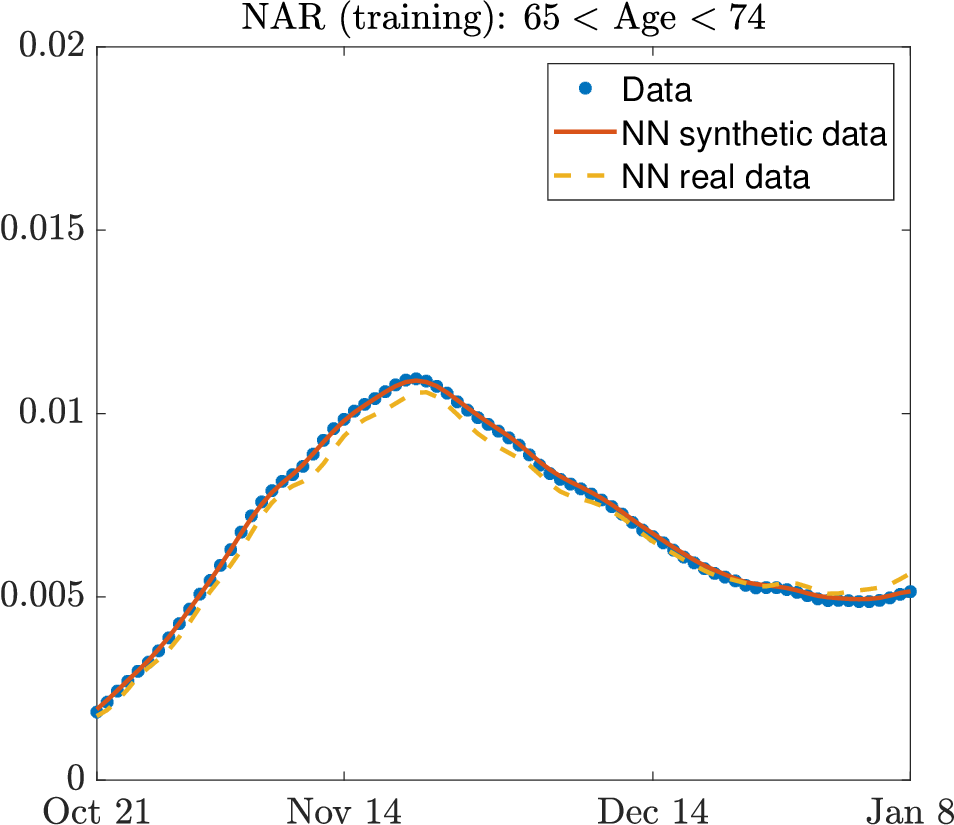}
	\includegraphics[width=0.328\linewidth]{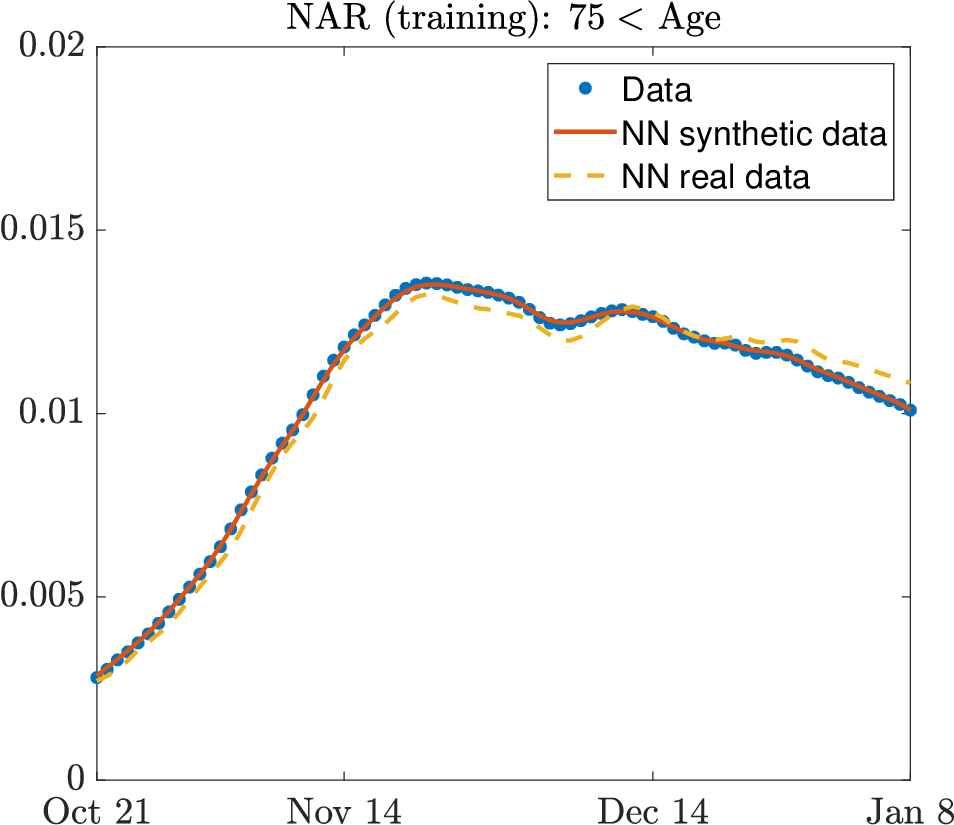}
	\caption{Nonlinear autoregressive network (training set) for the age-structured social SIAR model \eqref{eq:social_SIAR_ages}. Solution obtained by training a NAR network on both real and synthetic data, compared to the available data. Each image corresponds to a different age class.}
	\label{fig:NAR_ages_training}
\end{figure}

\begin{figure}[ht] 
	\centering
	\includegraphics[width=0.328\linewidth]{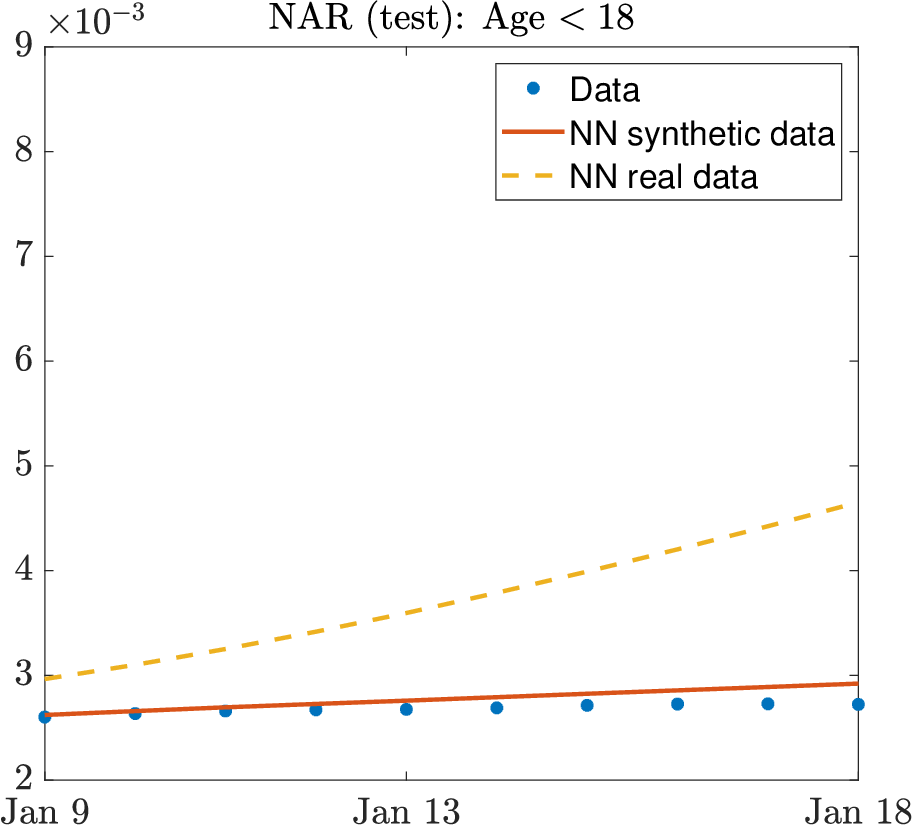}
	\includegraphics[width=0.328\linewidth]{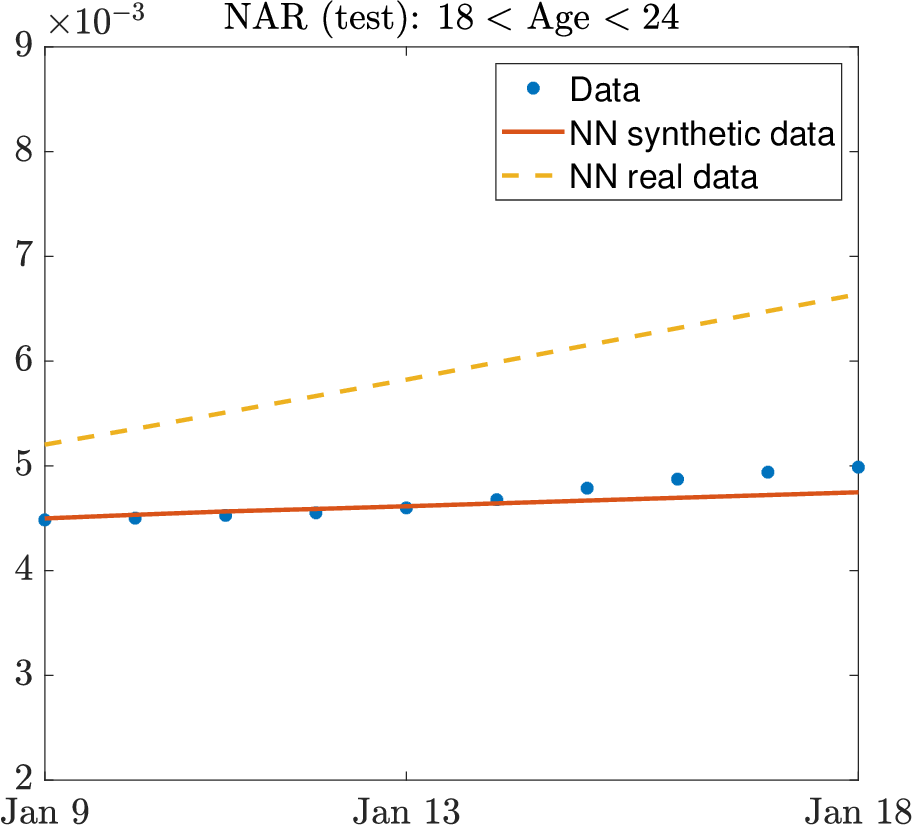}
	\includegraphics[width=0.328\linewidth]{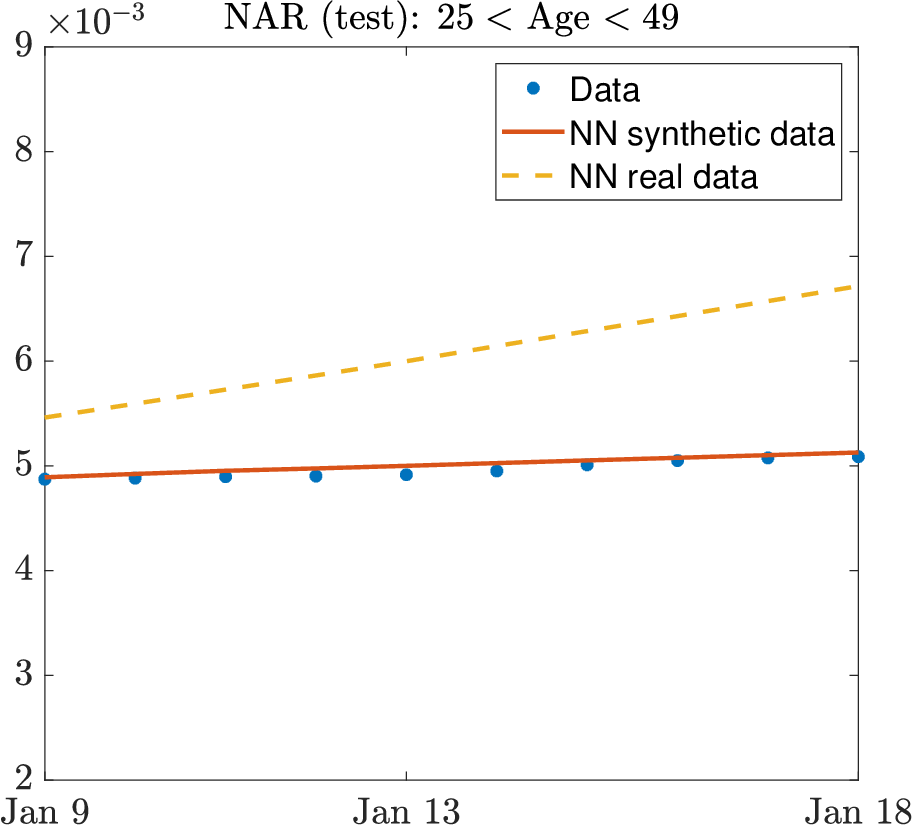}\\\vspace{0.2cm}
	\includegraphics[width=0.328\linewidth]{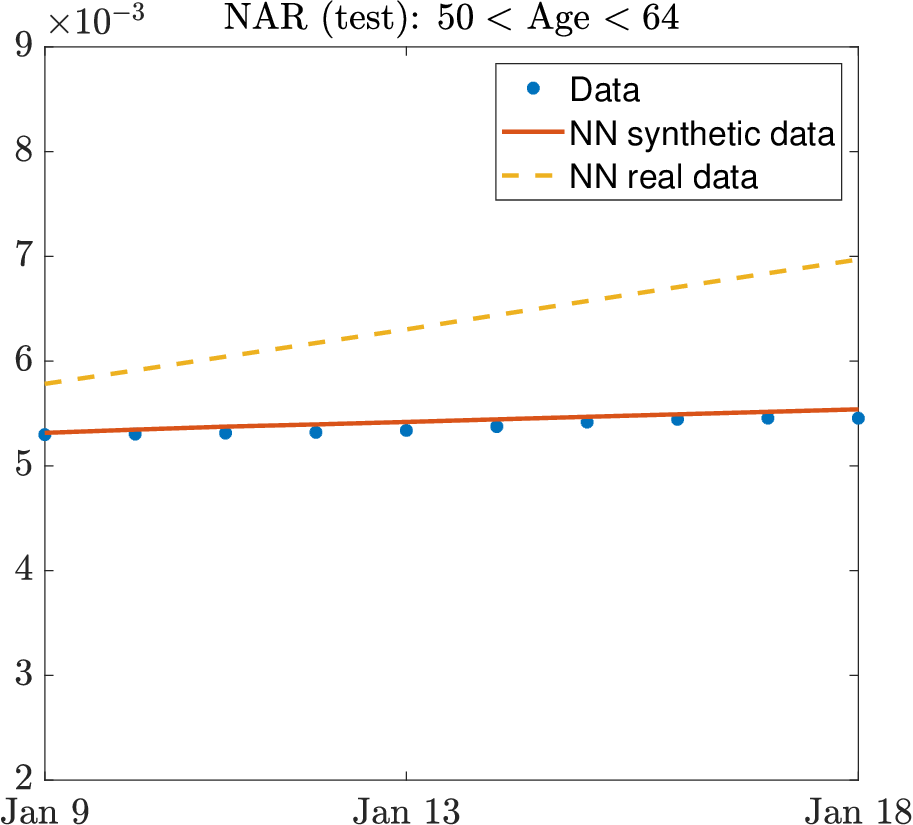}
	\includegraphics[width=0.328\linewidth]{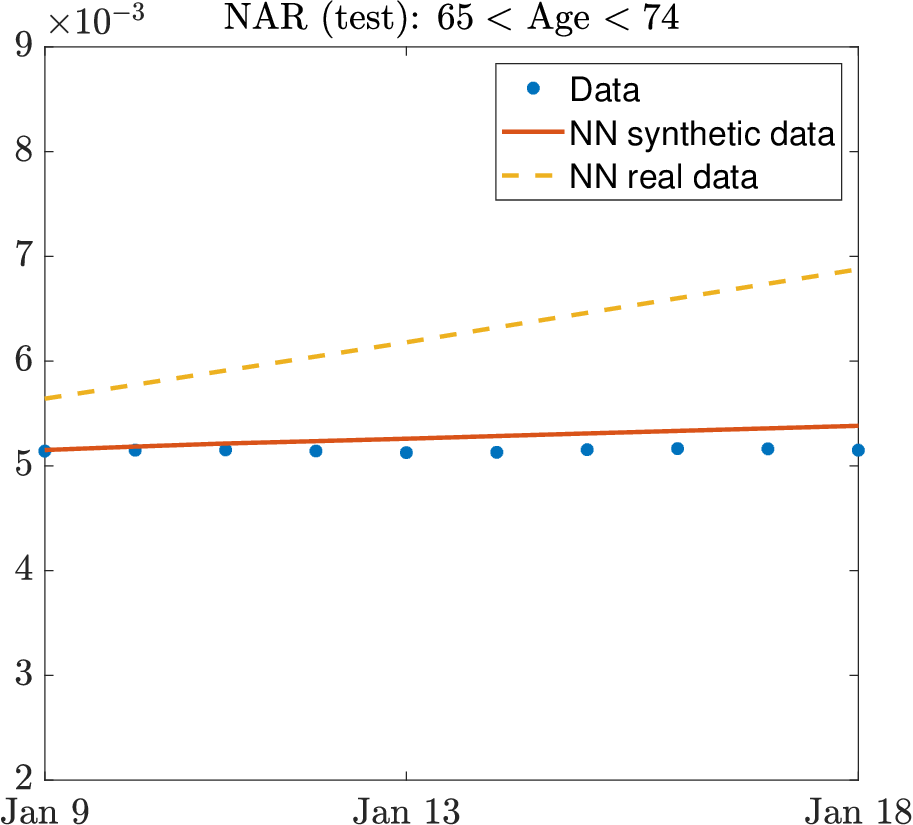}
	\includegraphics[width=0.329\linewidth]{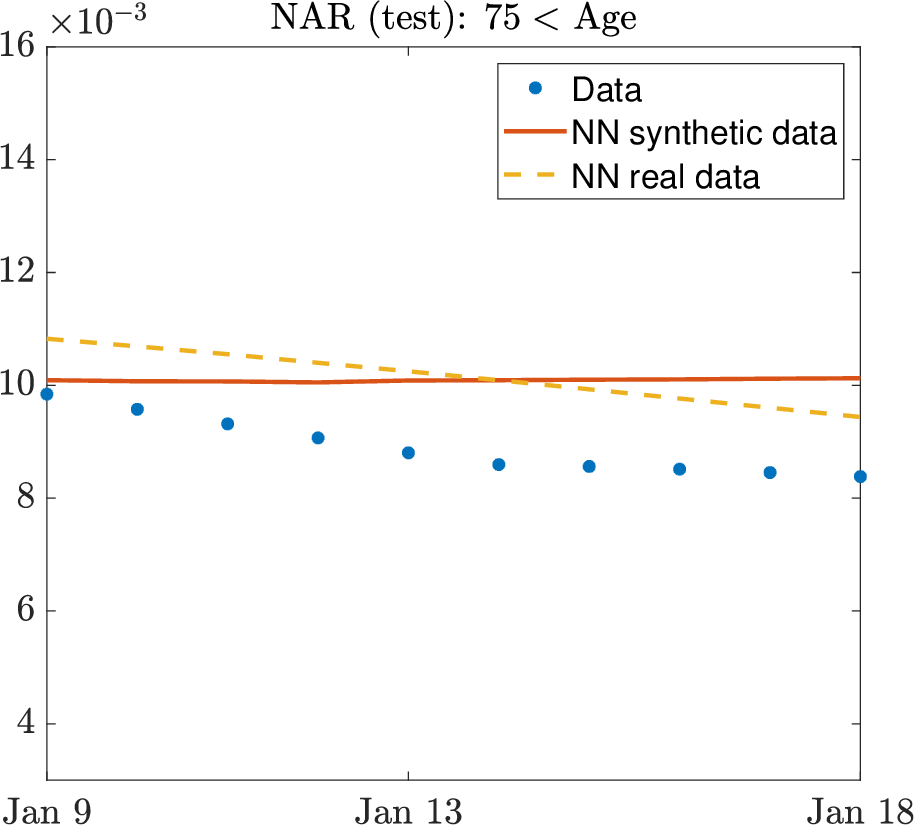}
	\caption{Nonlinear autoregressive network (test set) for the age-structured social SIAR model \eqref{eq:social_SIAR_ages}. Solution obtained by training a NAR network on both real and synthetic data, compared to the available data. Each image corresponds to a different age class.}
	\label{fig:NAR_ages_test}
\end{figure}

We now consider the NAR networks. We define two additional Feed-Forward network architectures  with width 32 and 3 hidden layers, using \texttt{ReLu} as activation function. One network is designed to approximate the solution of the social SIAR model \eqref{eq:social_SIAR}, and the other for the age-structured SIAR model \eqref{eq:social_SIAR_ages}.  As for the PINNs, both networks are  trained on real and synthetic datasets, assuming in the loss function \eqref{eq:loss}, $\omega_d = 1$ and $\omega_p=0$. In the case of real data, and for the non-age structured model, the network takes as input the number of infected individuals at times $t-d, \ldots, t-1$ within the training set and predicts the number of infected at time $t$, for any $t$ corresponding to real data. In the case of synthetic data, the input includes the number of infected at time $t-d,\ldots,t-1$ for any value of $z$, and the output is the number of infected for any value of $z$ at time $t$, for any $t$ corresponding to synthetic data. To select the lag parameter $d$, we performed a series of preliminary tests on the simpler model in \eqref{eq:social_SIAR}, evaluating the accuracy of the predicted solutions. Specifically, we let $d$ to vary in the range $d = 3, \ldots, 10$, and we measure the Root Square Mean Errors (RMSE) between the available data and the neural network predictions, defined as
\begin{equation}\label{eq:RMSE}
	RMSE(t) = \sqrt{\sum_{m=1}^{M} (\hat{I}(t,z_m) - I_i^{\text{NN}}(t,z_m;\theta_*))^2 w_m},
\end{equation}
where $\hat{I}(t,z_m)$ denotes the observed data and $I_i^{\text{NN}}(t,z_m;\theta_*)$ is the output of the neural network trained on real or synthetic data, and $w_m$ are the nodes associated to the uncertain parameters $z_m$. The results show that, for both synthetic and real data, changing the value of $d$ does not lead to a significant improvement in the accuracy of the neural network predictions, with mean RMSE being of the order of $10^{-4}$ for synthetic data and of the order of $10^{-3}$ for real data. Therefore, in all subsequent experiments, we set $d = 5$.  For the age-structured model, both the input and output include the infected number for all age classes. The networks are trained for 20000 epochs using Adam optimizer with a learning rate of $10^{-2}$.  To assess their performance, we first compute the solutions on the training set. Then, using a closed-loop strategy, where the network predictions are recursively used as inputs—we compute the solutions on the test set.
Figure \ref{fig:NAR_noAges} shows the comparison between the NAR network predictions and the available data in the case of the social SIAR model \eqref{eq:social_SIAR}. For synthetic data, we plot the mean solution with respect to the uncertainty, computed as in \eqref{eq:mean_sol}.
The image on the left corresponds to the training set, while on the right, the results on the test set are shown. The network trained on synthetic data yields more accurate solutions, demonstrating qualitative improved performance in both interpolation and prediction with respect to PINNs.
Figure \ref{fig:NAR_ages_training} presents the comparison between the neural network predictions and the available data for the age-structured social SIAR model \eqref{eq:social_SIAR_ages} computed on the training set while  Figure \ref{fig:NAR_ages_test} shows the corresponding results on the test set. Each plot corresponds to a different age class. In most of the cases, the network trained on synthetic data produces qualitatively more accurate results with respect to the one trained on real data and to PINNs, both in terms of interpolation and prediction. 

\subsection{Comparison between PINNs and NARs.} \label{sec:comparison}
\begin{figure}[ht] 
	\centering
	\includegraphics[width=0.328\linewidth]{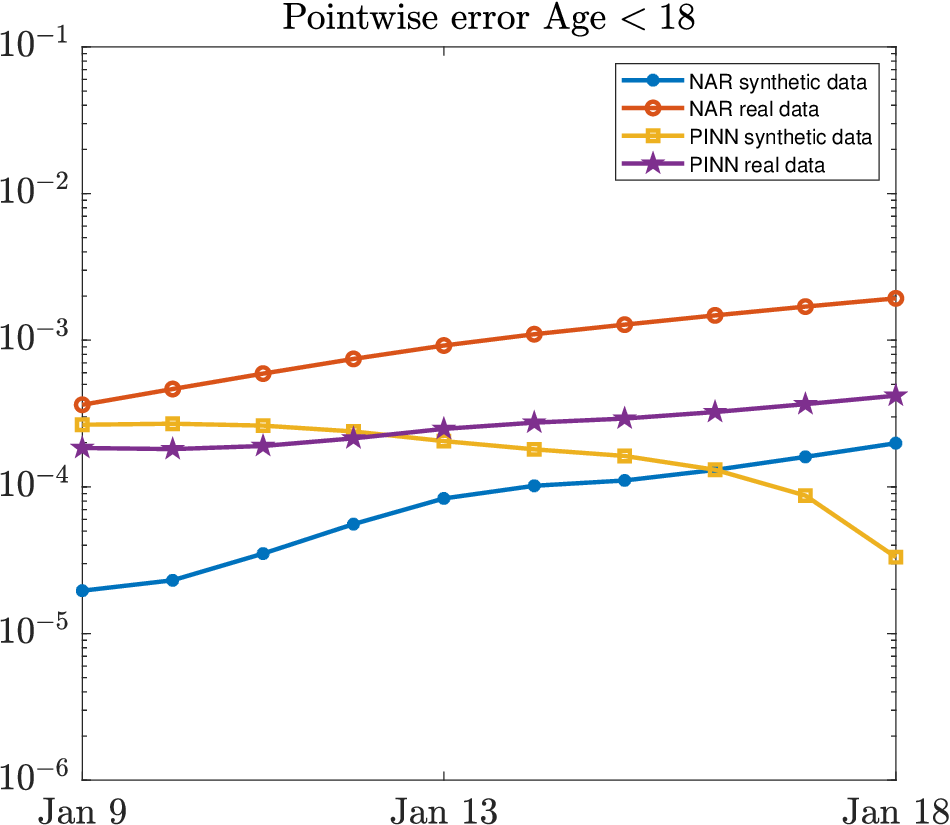}
	\includegraphics[width=0.328\linewidth]{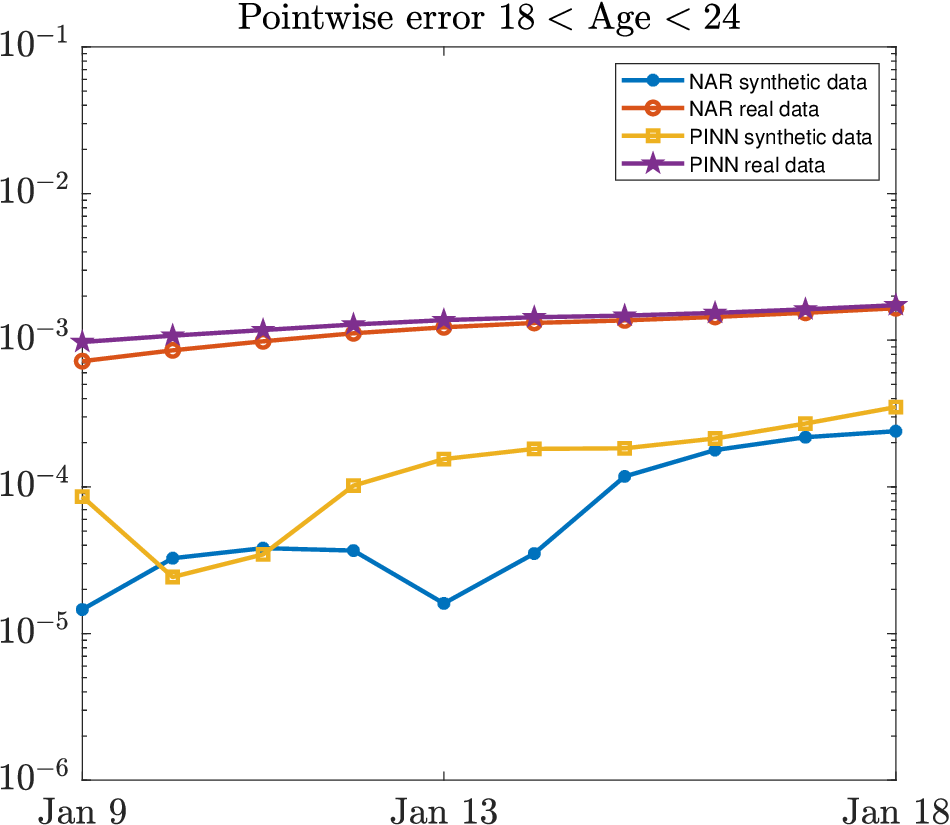}		\includegraphics[width=0.328\linewidth]{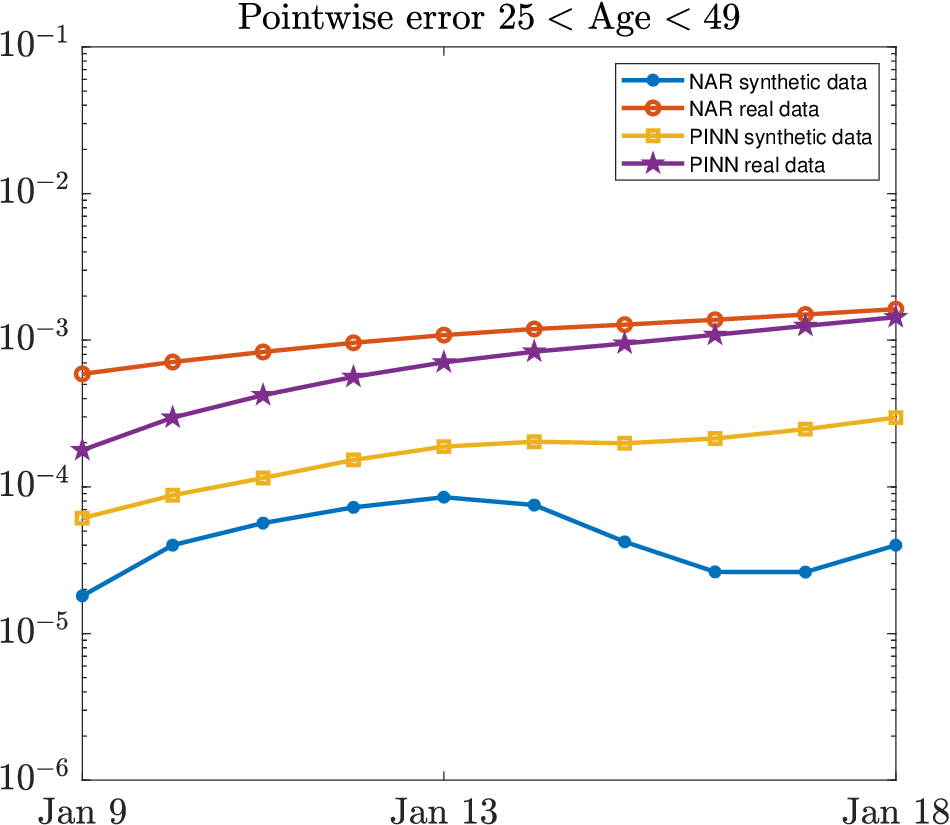}\\\vspace{0.2cm}
	\includegraphics[width=0.328\linewidth]{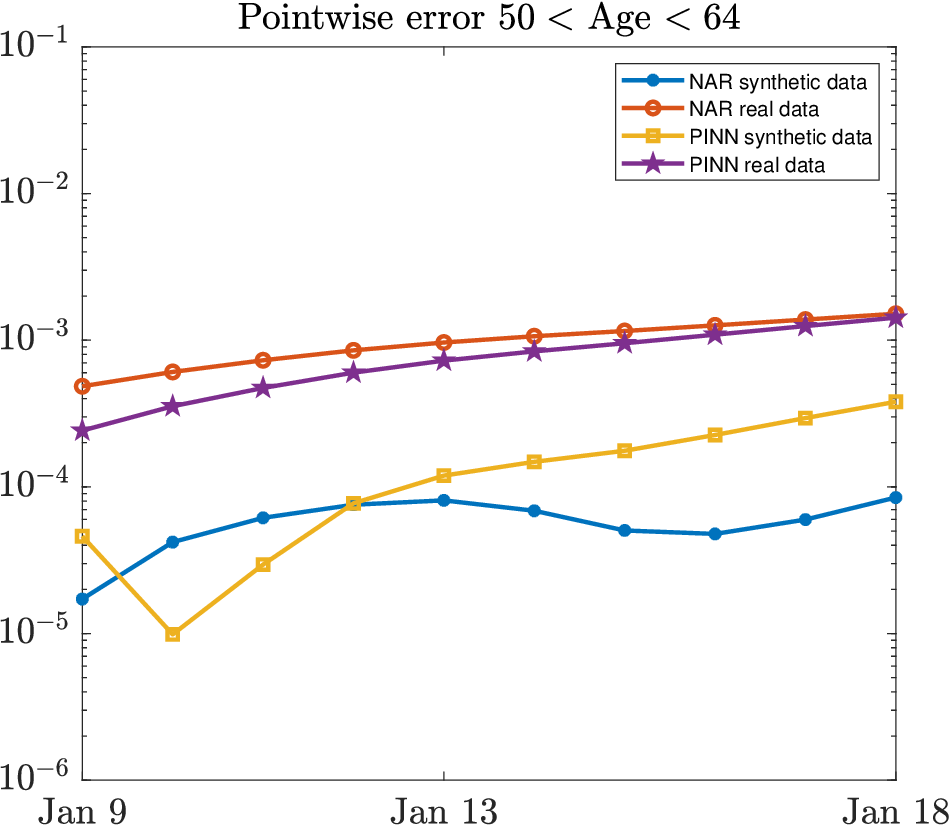}
	\includegraphics[width=0.328\linewidth]{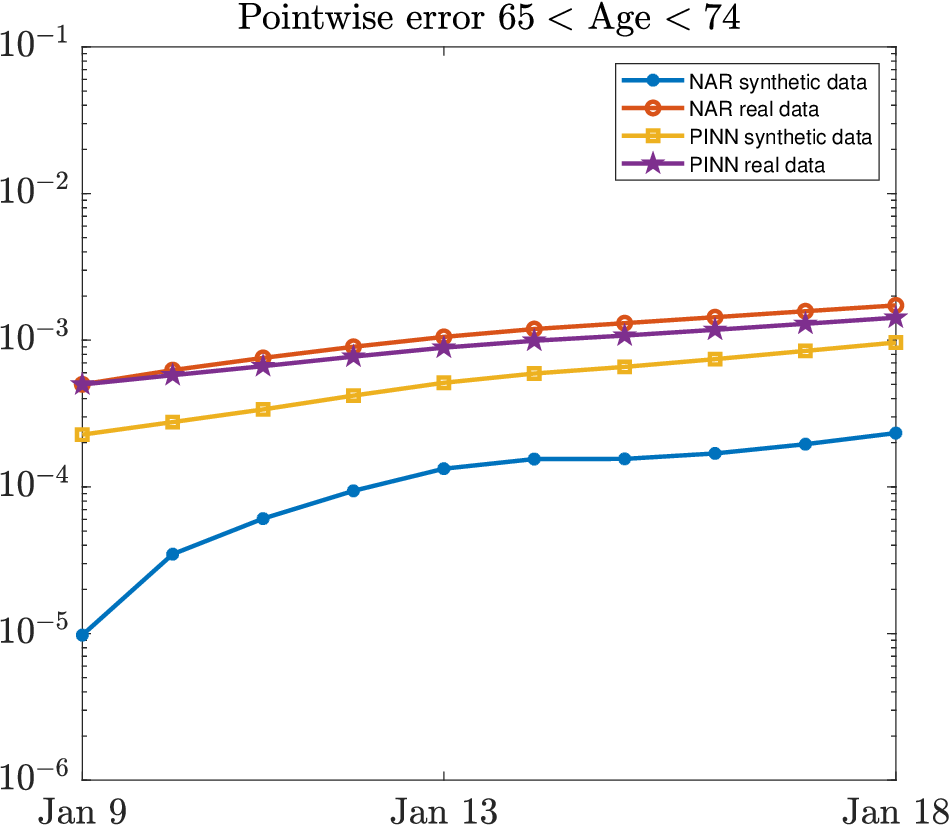}		\includegraphics[width=0.328\linewidth]{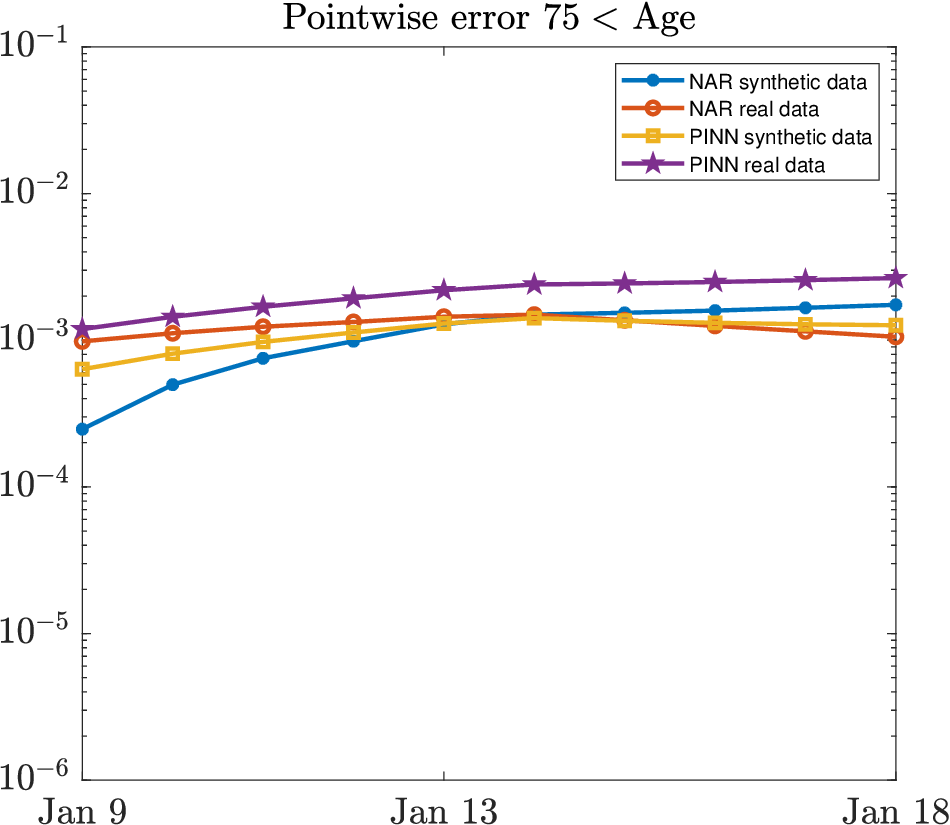}
	\caption{Comparison between NAR and PINN in terms of accuracy over the test set referred to the age-structured social SIAR model \eqref{eq:social_SIAR_ages}. Short term forecasting. The error between available data and neural network solutions in time is computed as in \eqref{eq:error}. Markers have been added just to denote different lines. }
	\label{fig:accuracy_ages}
\end{figure}

We now compare the performance of the NAR network and the PINN. We focus on both short term forecasting.  We first compare the trained networks in terms of computational cost. For comparison purposes, we assume to train PINNs and NAR networks, both on synthetic and real data, over 50000 epochs. Training a PINN involves significantly higher computational costs with respect to training NAR networks, particularly in the case of the age-structured model. Indeed, automatic differentiation is needed to enforce the physical constraints imposed by the system of ODEs. The results are confirmed by Table \ref{tab:costs} which reports the training time, in seconds, required to complete 50000 epochs referred to both models \eqref{eq:social_SIAR}-\eqref{eq:social_SIAR_ages}.  
\begin{table}[H]
	\begin{center}
		\caption{Training time in seconds required to complete 50000 epochs in the cases of NAR and PINN trained both on synthetic and real data. First row: social SIAR model \eqref{eq:social_SIAR}. Second row: age-structured social SIAR model \eqref{eq:social_SIAR_ages}. }
		\begin{tabular}{ccccc}
			\hline
			& NAR (synthetic) & NAR (real) & PINN (synthetic) & PINN (real)\\
			\hline
			Non-aged model  & 32s & 26s & 307s & 205s \\
			\hline
			Age-structured model   & 51s & 34s &468s & 234s \\
			\hline
		\end{tabular}
		\label{tab:costs}
	\end{center}
\end{table}   
\paragraph{Short term forecasting.}
In short-term forecasting, NAR networks offer a competitive alternative to physics-informed approaches, particularly when used in combination with data augmentation strategies. By directly learning the epidemic dynamics from available data, NAR networks provide accurate quantitative estimates that can be highly valuable for monitoring the progression of the pandemic and supporting decision-making aimed at mitigating its impact. Their data-driven nature also allows for faster training and computational efficiency, making them particularly suitable for timely analyses in rapidly evolving scenarios. The results are confirmed by 
Figure \ref{fig:accuracy_ages} which displays the error between the real data and the neural network solutions computed as \begin{equation}\label{eq:error}
	\mathcal{E}(x,t) = \bigl\lvert \hat{I}(x,t) - I_i^{\text{NN}}(x,t;\theta_*) \bigr\rvert,
\end{equation}
where $\hat{I}(x,t)$ denotes the observed data and $I_i^{\text{NN}}(x,t;\theta_*)$ is the output of the neural network trained on real or synthetic data. In the case of synthetic data, the network output corresponds to the mean solution with respect to the uncertainty, as in \eqref{eq:mean_sol}. 
Training either a PINN or a NAR network on synthetic data provides clear benefits in terms of accuracy. In particular, NAR networks trained on synthetic data proved to be effective in nearly all cases.
We present the results only for the age-structured model \eqref{eq:social_SIAR_ages}; however, similar findings hold for the simpler model \eqref{eq:social_SIAR} as well. These results are further supported by Table \ref{tab:accuracy}, which reports the maximum error \eqref{eq:error} over time between real data and the different neural network solutions. In the case of model \eqref{eq:social_SIAR} we integrate over $x\in \mathcal{A}$, being $\mathcal{A}$ the age-classes.
\begin{table}[H]
	\begin{center}	
		\caption{Accuracy of the neural network solutions w.r.t. the available data computed as the maximum of $\mathcal{E}(t)$ over time. Short term forecasting. First row: social SIAR model \eqref{eq:social_SIAR}. Second-Seventh rows: age-structured social SIAR model \eqref{eq:social_SIAR_ages}.} 
		\begin{tabular}{ccccc}
			\hline
			& NAR (synthetic) & NAR (real) & PINN (synthetic) & PINN (real)\\
			\hline
			Non-aged model & $2.1\times 10^{-4}$ & $1.2\times 10^{-3}$ & $10^{-3}$ &$1.2\times 10^{-3}$  \\
			\hline
			Age $<$ 18 & $2\times 10^{-4}$ & $1.9\times 10^{-3}$ &$3\times 10^{-4}$ & $4\times 10^{-4}$ \\
			\hline
			19$<$ Age $<$ 24 & $2\times  10^{-4}$ & $1.7\times 10^{-3}$ &$3\times 10^{-4}$ & $1.7\times 10^{-3}$ \\
			\hline
			25$<$ Age $<$ 49 & $10^{-4}$ & $1.6\times 10^{-3}$ &$3\times 10^{-4}$ & $1.4\times 10^{-3}$ \\
			\hline
			50$<$ Age $<$ 64 & $10^{-4}$ & $1.5\times 10^{-3}$ &$4\times 10^{-4}$ & $1.4\times 10^{-3}$ \\
			\hline
			65$<$ Age $<$ 74 & $2\times 10^{-4}$ & $1.7\times 10^{-3}$ &$  10^{-3}$ &$1.4\times 10^{-3}$ \\
			\hline
			75$<$Age & $1.7\times 10^{-3}$ & $1.5\times 10^{-3}$ &$1.4\times 10^{-3}$ &$2.7\times 10^{-3}$ \\
			\hline
		\end{tabular}
		\label{tab:accuracy}
	\end{center}
\end{table}
Let us finally mention that NAR networks achieve in general good accuracy after 20000 training epochs, whereas PINNs require at least 50000 epochs to produce good results. This is valid for both the social SIAR model and the one with age-structure.
\paragraph{Long term forecasting.}
In long-term forecasting, PINNs prove to be more accurate than NARs, by offering qualitative understanding of the epidemic dynamics. While their quantitative accuracy may be lower than that of NAR networks in the short term, PINNs are particularly useful for exploring strategies to mitigate epidemic peaks and understand broader trends in disease progression. By embedding the underlying physical and epidemiological laws into the learning process, they can capture the overall behavior of the system over extended periods, helping to inform public health interventions and long-term planning. Since the main goal of this test is to illustrate the capability of PINNs to reproduce the epidemic peak, we restrict our analysis to the simpler model in \eqref{eq:social_SIAR}. Capturing the additional complexity introduced by age-structured dynamics would likely require a significantly greater training effort from the network. Furthermore, as we have already shown that the augmented data strategy improves the accuracy of neural network solutions, we train both PINNs and NARs exclusively on augmented datasets, without reporting results obtained from training on real data.

\begin{figure}[h]
	\centering
	\includegraphics[width=0.426\linewidth]{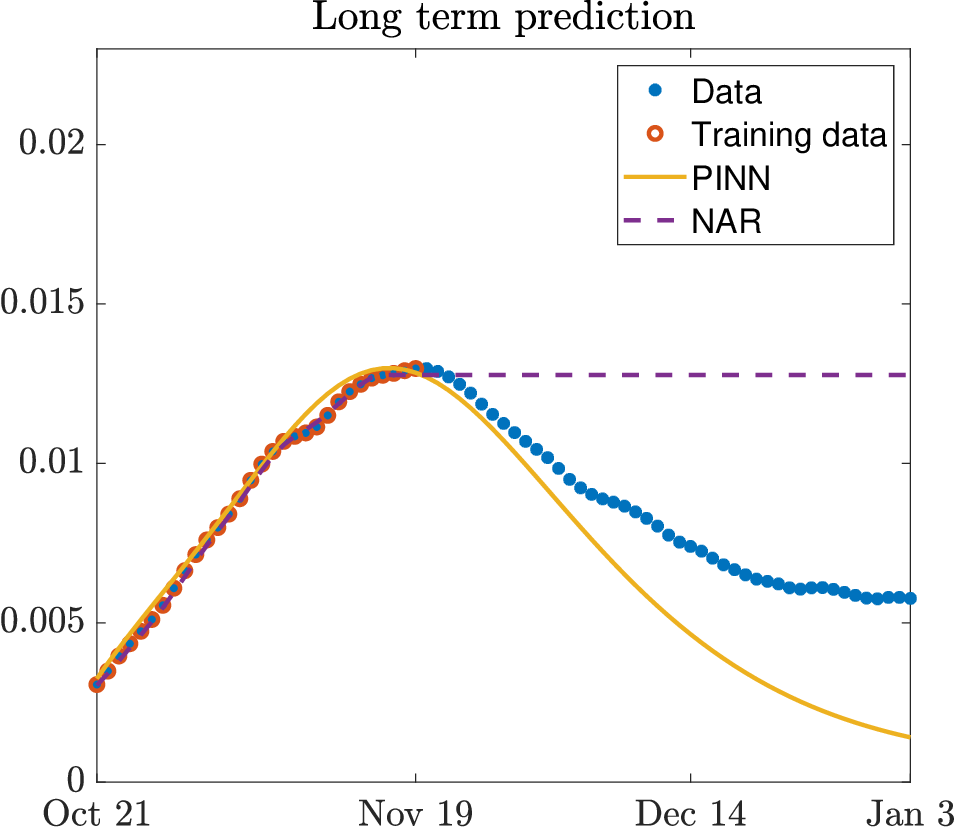}
	\includegraphics[width=0.415\linewidth]{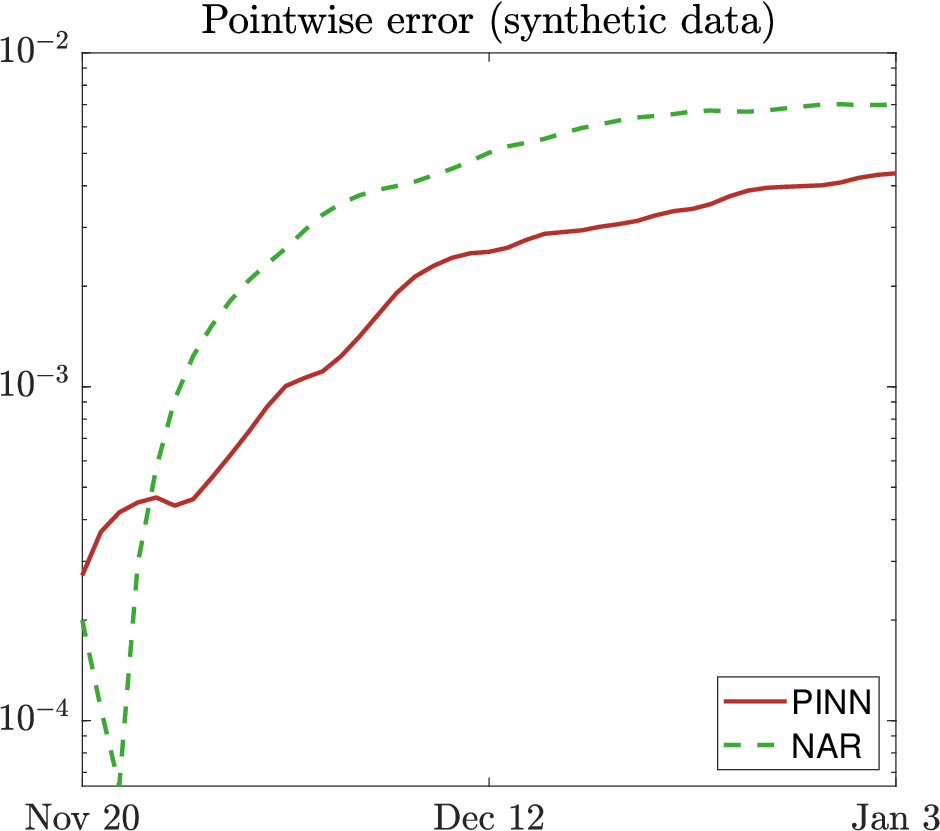}
	\caption{Long term forecasting for the  social SIAR model \eqref{eq:social_SIAR}. On the left, the solution obtained by training a PINN and a NAR network on synthetic data, compared to the available data. On the right, the pointwise error over the test set computed as in \eqref{eq:error} }
	\label{fig:accuracy_long_time}
\end{figure}

Figure \ref{fig:accuracy_long_time} (left) compares NAR and PINN networks trained on synthetic data for long-term forecasting. PINNs achieve higher accuracy than NARs, successfully reproducing the epidemic peak and capturing the qualitative behavior of the solution. These findings are further supported by the error plot, computed as in \eqref{eq:error} on the test set, shown in the right panel of Figure \ref{fig:accuracy_long_time}.

All the numerical experiments have been run on a Desktop Computer equipped with an Intel(R) Core(TM) i7-8700 CPU processor and 32GB RAM.

\subsection{Effectiveness of the data augmentation strategy} 	
\begin{figure}[h]
	\centering
	\includegraphics[width=0.424\linewidth]{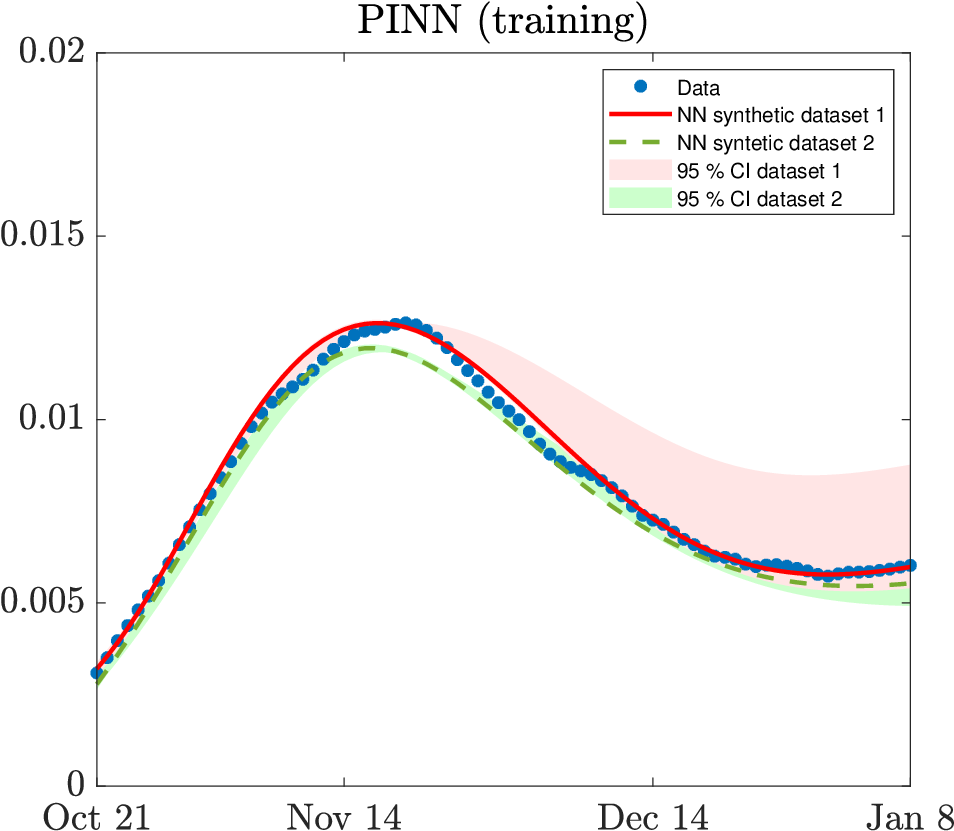}
	\includegraphics[width=0.41\linewidth]{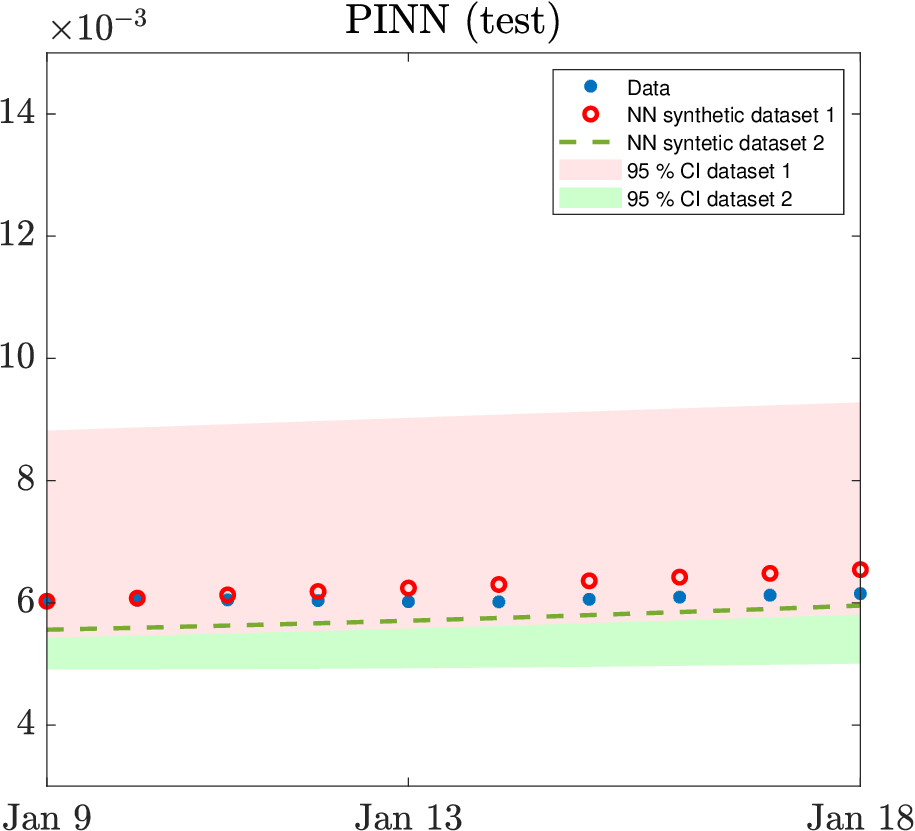}\\
	\includegraphics[width=0.424\linewidth]{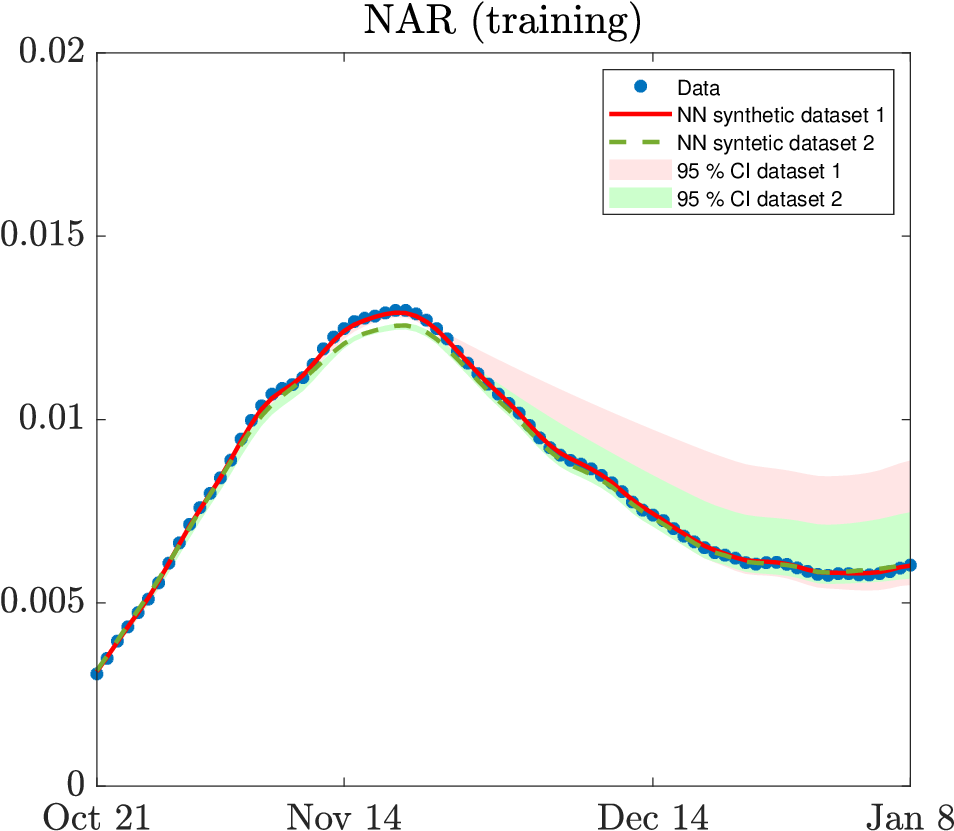}
	\includegraphics[width=0.41\linewidth]{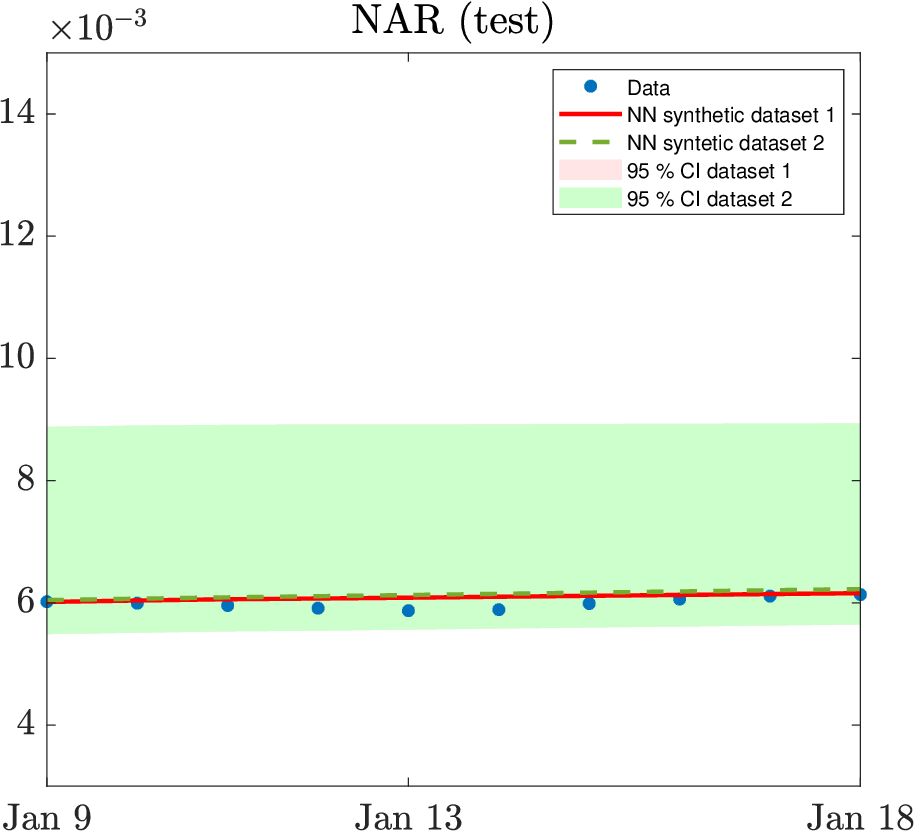}
	\caption{ NAR network and PINN for the social SIAR model \eqref{eq:social_SIAR}. Comparison between two different datasets of uncertain parameters computed as in \eqref{eq:gamma} and following the procedure described in Section \ref{sec:params_est}. First row: PINN. Second row: NAR networks. On the left, the solution computed on the training set. On the right, the solution computed on the test set.  Dataset one and two correspond to $z_1 \sim \mathrm{Beta}(1.95,4.95)$ and  $z_2 \sim \mathrm{Beta}(2,6)$ , respectively.      }
	\label{fig:PINN_NAR_bande}
\end{figure}

To conclude our study and assess the effectiveness of the proposed data augmentation strategy, we perform an additional test focusing on short-term forecasting, referring for simplicity to the model introduced in \eqref{eq:social_SIAR}. We assume the recovery parameters to follow beta distributions defined as
\begin{equation}\label{eq:gamma_distr} 
	\begin{split}		
		\gamma_I^1 = 0.04 + 0.05 z_1, \qquad \gamma_I^2 = 0.04 + 0.05 z_2,
	\end{split}
\end{equation} 	
where $z_1 \sim \mathrm{Beta}(1.95,4.95)$ and $z_2 \sim \mathrm{Beta}(2,6)$. Following the procedure described in Section~\ref{sec:params_est}, we set $\gamma_A^i = 2\gamma_I^i$ for $i=1,2$, and determine the transmission parameter $\beta$, the unknown number of asymptomatic $\xi(\cdot)$ and the macroscopic incidence rates $H_J(\cdot)$, for $J\in\{S,I,A\}$ by solving the corresponding minimization problems. 
As in the previous experiments, we generate synthetic data by solving the ODE system in \eqref{eq:social_SIAR} and evaluate the performance of both the PINN and NAR networks on this datasets. Specifically, we select the period from October 21st to January 8th as the training set, and the interval from January 9th to 18th as the test set. Figure \ref{fig:PINN_NAR_bande} illustrates the solutions computed by the two networks trained on the respective datasets for the different uncertain parameters. Specifically, dataset one and two correspond to $z_1 \sim \mathrm{Beta}(1.95,4.95)$ and  $z_2 \sim \mathrm{Beta}(2,6)$ , respectively.   The plots display the mean predicted solutions, with the confidence intervals represented by the shaded area. The left panel corresponds to the training phase, while the right panel to the testing phase. In the first row the plots for the NAR network and in the second row for the PINN. In terms of mean solution, the results are comparable. 
Finally, we compute the accuracy of the solution on the test set according to different error metrics.
Specifically, we consider the Root Mean Square Error (RMSE) computed as in \eqref{eq:RMSE}, as well as the Mean Absolute Error (MAE) and the coefficient of determination ($R^2$),
\begin{equation}\label{eq:errors}
	\begin{aligned}
		&MAE(t) = \frac{1}{N_c} \sum_{m=1}^{M} \big|\hat{I}(t,z_m) - I_i^{\text{NN}}(t,z_m;\theta_*)\big|w_m, \\
		&R^2(t) = 1 - \frac{\sum_{m=1}^{M} \big(\hat{I}(t,z_m) - I_i^{\text{NN}}(t,z_m;\theta_*)\big)^2 w_m}{M\hat{\sigma}^2},
	\end{aligned}
\end{equation}
where $\hat{I}(t,z_m)$ denotes the observed data, $I_i^{\text{NN}}(t,z_m;\theta_*)$ is the output of the neural network trained on real or synthetic data, $N_c$ is the number of considered time points, and $\hat{\sigma}^2$ is the variance of the observed data with respect to the uncertainty parameter $z_m$, and $w_m$ are the nodes associated with $z_m$, for $m=1,\ldots,M$. Table~\ref{tab:errors_2} reports the mean value with respect to time $t$ of the quantities computed in \eqref{eq:RMSE}-\eqref{eq:errors}. The coefficient of determination $R^2$ indicates that the model provides an excellent fit to the data, a conclusion further supported by the other error metrics. Moreover, the neural networks demonstrate robustness, as the errors remain comparable across different uncertainty levels. Once again, the results confirm that, for short-term predictions, the NAR network yields more accurate approximations than the PINN.
\begin{table}[H]
	\begin{center}	
		\caption{ Accuracy of the neural network solutions to the social SIAR model \eqref{eq:social_SIAR} w.r.t. the available data computed as in \eqref{eq:RMSE}-\eqref{eq:errors}. Short term forecasting.} 
		\begin{tabular}{ccccc}
			\hline
			& NAR (dataset 1)  & PINN (dataset 1)& NAR (dataset 2)  & PINN (dataset 2)\\
			\hline
			$RMSE$ &$2\times 10^{-4}$ & $6.2 \times 10^{-4}$ & $10^{-4}$ & $6\times 10^{-4}$ \\
			\hline
			$MAE$ & $1.8\times 10^{-4}$ & $10^{-3}$ & $2.2\times 10^{-4}$ & $7\times 10^{-4}$\\
			\hline 
			$R^2$ & $0.99$ &$0.99$ &$0.99$ &$0.99$\\
			\hline
		\end{tabular}
		\label{tab:errors_2}
	\end{center}
\end{table}

\section{Conclusion}\label{sec:conclusion} 
In this work, we explored a data augmentation strategy aimed at enhancing the performance of neural networks for epidemic modeling, focusing on both interpolation and prediction tasks. By combining real-world data with synthetic data generated through simulations of socially structured SIAR models that incorporate parameter uncertainty and age dependence, we demonstrated the ability of neural networks to effectively capture complex epidemic dynamics when additional data from suitable models are included into the training. In addition to  PINNs, we considered an emerging class of models known as NAR networks. These networks, which use past time-step solutions to predict future dynamics, have shown particular effectiveness in capturing temporal dependencies, especially in the short time forecasting. Unlike PINNs, which require embedding the governing equations directly into the training process, NAR networks adopt a purely data-driven approach that has proven to be both more accurate and significantly less computational expensive. A key finding of this study is that incorporating uncertainty into model parameters not only enables the generation of more realistic synthetic datasets but also improves the predictive accuracy of neural networks. 

Future research will explore alternative neural network architectures within the class of recurrent models, such as Long Short-Term Memory (LSTM) networks, which may offer further advantages for time-series forecasting. Moreover, we plan to extend our framework to spatially dependent epidemic models under uncertainty. This direction is motivated by the fact that disease transmission is often spatially heterogeneous. Urban areas, for example, typically exhibit faster spread than rural regions, and spatially targeted interventions such as localized lockdowns or travel restrictions introduce additional complexity. Incorporating spatial structure and uncertainty will enable a more realistic and localized understanding of epidemic dynamics, supporting more effective and adaptive containment strategies.

\section*{Acknowledgments}
This work has been written within the activities of GNCS and GNFM groups of INdAM
(Italian National Institute of High Mathematics). G.D. has been partially funded by the European Union — NextGenerationEU, MUR–PRIN 2022 PNRR Project No. P2022JC95T “Datadriven discovery and control of multi-scale interacting artificial agent systems”. G.D. and F.F.
thank the Italian Ministry of University and Research (MUR) through the PRIN 2020 project
(No. 2020JLWP23) “Integrated Mathematical Approaches to Socio–Epidemiological Dynamics”. L.P. has been partially funded by the European Union– NextGenerationEU under the
program “Future Artificial Intelligence– FAIR” (code PE0000013), MUR PNRR, Project “Advanced MATHematical methods for Artificial Intelligence– MATH4AI”. L.P. acknowledges the
support by the Royal Society under the Wolfson Fellowship “Uncertainty quantification, datadriven simulations and learning of multiscale complex systems governed by PDEs” and by
MIUR-PRIN 2022 Project (No. 2022KKJP4X), “Advanced numerical methods for time dependent parametric partial differential equations with applications”. The partial support by
ICSC – Centro Nazionale di Ricerca in High Performance Computing, Big Data and Quantum
Computing, funded by European Union – NextGenerationEU is also acknowledged.

\bibliographystyle{abbrv}
\bibliography{biblio}
\end{document}